\documentclass[a4paper,12pt,reqno]{amsart}
\usepackage{amssymb}

\newtheorem{lemma}{Lemma}
\newtheorem{proposition}[lemma]{Proposition}
\newtheorem{theorem}[lemma]{Theorem}
\newtheorem{corollary}[lemma]{Corollary}

\newtheorem{comm}[lemma]{Comment}
\newtheorem{remm}[lemma]{Remark}
\newtheorem{exam}[lemma]{Example}

\newenvironment{example}{\begin{exam} \fontsize{10pt}{7} \fontshape{n}\selectfont}
{\end{exam}}
\newenvironment{vect}{\left( \!\!\! \begin{array}{c}}{\end{array} \!\!\!\right)}
\newenvironment{matr2}{\left( \!\!\! \begin{array}{cc}}{\end{array} \!\!\!\right)}

\newcommand{\alp}{\alpha}
\newcommand{\bet}{\beta}
\newcommand{\gam}{\gamma}

\newcommand{\lam}{\lambda}
\newcommand{\ome}{\omega}
\newcommand{\del}{\delta}
\newcommand{\eps}{\varepsilon}
\newcommand{\tet}{\vartheta}

\newcommand{\Ome}{\Omega}
\newcommand{\Lam}{\Lambda}

\renewcommand{\Re}{\mathrm{Re}}
\renewcommand{\Im}{\mathrm{Im}}

\newcommand{\R}{\mathbb{R}}
\newcommand{\C}{\mathbb{C}}
\newcommand{\Q}{\mathbb{Q}}
\newcommand{\Z}{\mathbb{Z}}
\newcommand{\Lp}[1][2]{{\rm L}^{#1}}
\newcommand{\So}[1][1,2]{{\rm W}^{#1}}

\newcommand{\Schrodinger}{Schr\"o\-din\-ger }

\newcommand{\ud}{{\rm d}}
\newcommand{\e}{{\rm e}}

\newcommand{\hull}{\mathrm{\hull}}

\newcommand{\spec}[2][\!]{\mathrm{Spec}_{#1} \left( #2 \right) \, }

\newcommand{\nume}[1]{\mathrm{Num} \left( #1 \right) \, }

\newcommand{\dist}[1]{\mathrm{dist} \left( #1 \right) \, }

\newcommand{\dotp}[1]{\left< #1 \right>}
\newcommand{\odtp}[1]{\left< #1 \right>}
\newcommand{\rang}[1]{\mathrm{Ran} \left( #1 \right) \, }

\newcommand{\dom}{\mathrm{Dom}}

\newcommand{\gap}{\vspace{.2in}}
\newcommand{\Proof}[1][\!\!]{\textbf{\underline{Proof}.}\, #1}
\newcommand{\fin}{\hspace*{\fill} $\, \blacksquare$}
\newcommand{\find}{\hspace{.1in} \blacksquare}

\fontshape{n}

\title[Spectral behaviour of a simple n-s-a operator]{SPECTRAL
    BEHAVIOUR OF \break
   A SIMPLE NON-SELF-ADJOINT OPERATOR}
\author[Lyonell S. Boulton]{}

\date{$6^{th}$ March 2001}

\subjclass{34L05, 47E05, 34L16.}

\keywords{Spectral theory of non-self-adjoint operators,
differential operators, non-real eigenvalues.}

\thanks{$\dag$ The author carried out this research as a PhD student at King's
College London sponsored by ``Fundaci\'on Gran Mariscal de
Ayacucho'', Venezuela, grant E-211-1357-1997-1.}

\begin{document}
\

\vskip -1in {\em \
\parbox[t]{1.1\linewidth}{\small
\rightline{KCL-MTH-01-07} } }

\vskip -1.5em \null

\vskip 1.5em

\maketitle

\centerline{\scshape   Lyonell S. Boulton$^\dag$}

\medskip

{\footnotesize \centerline{ Department of Mathematics, King's
College London } \centerline{The Strand, London WC2R 2LS, U.K. }
\centerline{email: lboulton@mth.kcl.ac.uk}}

\medskip

\begin{abstract}
We investigate the spectrum of a typical non-self-adjoint
differential operator $AD=-d^2/dx^2\otimes A$ acting on \linebreak
$\Lp(0,1) \otimes \mathbb{C}^2$, where $A$ is a $2\times 2$
constant matrix. We impose Dirichlet and Neumann boundary
conditions in the first and second coordinate respectively at both
ends of $[0,1]\subset \R$. For $A\in \R^{2\times 2}$ we explore in
detail the connection between the entries of $A$ and the spectrum
of $AD$, we find necessary conditions to ensure similarity to a
self-adjoint operator and give numerical evidence that suggests a
non-trivial spectral evolution.
\end{abstract}

\section{Introduction} \label{s1}

In this paper we investigate spectral properties of the linear
operator $AD$ acting on $\Lp(0,1) \otimes \C^2$ where $A$ is a $2
\times 2$ constant matrix and $D$ denotes the ordinary
differential operator
\[
D \begin{vect} \phi \\ \gam \end{vect} := -\begin{vect}\phi'' \\
\gam''
\end{vect}, \qquad \qquad \begin{aligned} \phi(0) & =\phi(1)=0 \\
\gam'(0)& =\gam'(1)=0. \end{aligned}
\]
The apparently simple combination of Dirichlet and Neumann
\linebreak boundary conditions allows self-adjointness if, an only
if, $A$ is real and diagonal. If $A$ is non-diagonal and
upper-triangular the numerical range of $AD$ is a large sector of
$\C$. Otherwise it is the whole of $\C$ preventing us from
applying the theory of sectorial sesquilinear forms in a
straightforward manner.

Our main goal is to explore the connection between the entries of
the matrix $A$ and the location of the spectrum of $AD$ in the
complex plane. In \cite{S1} R.~F.~Streater considers the
particular case
\[
    A=\begin{matr2} 1 & \gam \\ 1/2\gam & 1 \end{matr2}, \qquad
    \qquad \gam>0,
\]
in order to find necessary conditions for the stability of small
perturbations about the stationary solution of certain non-linear
system of parabolic equations. Streater's system represents a
thermodynamical model for hot fluid in one dimension and the
localization of the spectrum is achieved by constructing a
non-unitary transformation that makes $AD$ similar to a
non-negative self-adjoint operator, hence the spectrum of $AD$ is
real and non-negative. This similarity transformation does not
work for other matrices and a slight modification of the entries
of $A$ can destroy reality of the spectrum (cf. sections
\ref{s6}-\ref{s7}) so the general case should be attacked by other
methods.

Although this paper mainly concerns $A\in \R^{2\times 2}$, the
results of sections \ref{s2}-\ref{s5} refer to any complex
$2\times 2$ matrix. The core results are to be found in section
\ref{s6} where we present an exhaustive description of the
spectrum of $AD$ in terms of the entries of $A$. Among various
other unexpected conclusions, the following three epitomize the
complexity of the problem to be considered:
\begin{itemize}
\item[a)] When $A$ is triangular and non-diagonalizable, $AD$ is
not similar to a self-adjoint operator but the spectrum of $AD$ is
real (theorem \ref{t11}).
\item[b)] The spectrum of $AD$ can be non-real even when both
 eigenvalues of $A$ are positive and equal
(theorem \ref{t32}).
\item[c)] There is a continuous family of matrices $A$ whose
eigenvalues do not intersect the real line but such that the
spectrum of $AD$ is real (theorem \ref{t31}).
\end{itemize}
The last two assertions show that the spectra of $A$, $D$ and $AD$
are typically unrelated.

The crucial idea in section \ref{s6} is to reduce the
four-parameter problem of localizing the spectrum of $AD$ in terms
of the entries of $A$, to five two-parameter cases and describe
separately each of these cases. Sections \ref{s2}-\ref{s5} are
devoted to describing the various properties of $AD$ we will use
in section \ref{s6}, whereas section \ref{s7} is devoted to
numerical computations which illustrate some of the results
reported. In section \ref{s2} we find the boundary conditions
associated to the adjoint of $AD$ and compute the numerical range
of $AD$. In section \ref{s3} we show that the resolvent of $AD$ is
compact for all non-singular $A$. In Section \ref{s4} we explore
the stability of the spectrum of $AD$ in the sense of \cite{D2}
and \cite{T1}, and provide estimates which allow us to enclose the
spectrum of $AD$ in angular regions when $A$ is subject to various
constraints. In section \ref{s5} we use standard ODE methods to
compute the transcendental function of the spectral problem
associated to $AD$.

\gap

\section{Definitions and notation} \label{s2}

Let $K$ be a linear operator whose domain is denoted by $\dom
(K)$. Throughout this paper $\spec{K}$ stands for the spectrum of
$K$ and the numerical range of $K$ is defined to be
\[
   \nume{K}:= \{ \dotp{Kf,f} : f \in \dom (K), \, \|f\|=1 \}.
\]
We recall that the numerical range of any linear operator is
convex and that if $\spec{K}\not=\emptyset$, then
\[
   \spec{K} \subset \overline{\nume{K}}.
\]

If $K=K^\ast$ and $\spec{K}\subset (0,\infty)$, we will say that
$K$ is positive and write $K>0$. If $K=K^\ast$ and
$\spec{K}\subset [0,\infty)$, we will say that $K$ is non-negative
and write $K\geq 0$.

Below and elsewhere $|v|$ denotes the norm of a vector $v\in
\C^2$. The norm of any
\[
f\equiv \begin{vect} \phi\\ \gam \end{vect} \in\Lp(0,1)\otimes
\C^2
\]
is the standard Hilbert tensor product norm
\[
\|f\|^2=\dotp{f,f}=\int_0^1 |f(x)|^2 \ud x = \int_0^1
\left(|\phi(x)|^2+|\gam(x)|^2\right) \ud x.
\]

Unless explicitly stated, we denote
\[
   A:= \begin{matr2} a & b \\ c & d \end{matr2}.
\]
The complex numbers $a_+,a_-$ denote the eigenvalues of $A$ and
the non-zero $\C^2$ vectors $v_+,v_-$ denote the eigenvectors
\[
Av_\pm=a_\pm v_\pm.
\]
If $a_+$ and $a_-$ are real and different, we adopt the convention
$a_-< a_+$.  Notice that the $v_\pm$ are not necessarily
orthogonal.

Let $\So[2,2]$ be the Sobolev space of all $f\in \Lp(0,1)\otimes
\C^2$, such that the generalized derivative $f''\in
\Lp(0,1)\otimes \C^2$. We define rigorously the domain of $AD$ as
\[
\dom(D)=\{ f \in \So[2,2] \, : \, \phi(0)=\phi(1)=0,\,
\gam'(0)=\gam'(1)=0\}.
\]
If $A$ is invertible, it is standard to show that $AD$ is always a
closed densely defined linear operator acting on $\Lp(0,1)\otimes
\C^2$.

\begin{lemma} \label{t27}
If $A$ is singular, then $AD$ is not closed in the domain
$\dom(D)$.
\end{lemma}
\Proof Let $v\in \C^2$ be a non vanishing vector such that $Av=0$
and let $f(x):=vx \in \Lp (0,1) \otimes \C^2$. Clearly $f \not \in
\dom(D)$. Let $\phi_n$ be a sequence of smooth functions whose
support is compact in $(0,1)$ and such that $\phi_n(x) \to x$ in
$\Lp(0,1)$. Then $\phi_n v \in \dom(D)$ and $\phi_n v \to f$. Also
\[
   AD\phi_n(x)v = -\phi_n''(x) Av =0,
\]
so that $AD(\phi_n v)$ is a convergent sequence in $\Lp(0,1)
\otimes \C^2$. We complete the proof by noticing that if $AD$ was
closed, then we would have $f\in \dom(D)$. \fin

For the rest of this section and in sections \ref{s3}-\ref{s5} we
will assume without further mention that $A$ is non-singular. In
section \ref{s6} we will consider again singular $A$.

\gap

In order to show that $AD$ is in general non-self-adjoint, let us
compute the adjoint $(AD)^{\ast}$. Let
\[
   P:= \begin{matr2} 1 & 0 \\ 0 & 0 \end{matr2}.
\]
Then the boundary conditions for $D$ can be rewritten as
\[
   Pf(0)=Pf(1)=0, \qquad (I-P)f'(0)=(I-P)f'(1)=0.
\]

\begin{lemma} \label{t18}
The adjoint of $AD$ is
\[
   (AD)^\ast f = -A^\ast f'',
\]
for $f\in \So[2,2]$ subject to the boundary conditions
\begin{equation} \label{e1}
   \begin{gathered}
   \hat{P}f(0)=\hat{P}f(1)=0 \\
   (I-\hat{P})f'(0)=(I-\hat{P})f'(1)=0
   \end{gathered}
\end{equation}
where $\hat{P}=\hat{P}^2$ is the rank one projection such that
\begin{equation} \label{e25}
   \begin{gathered}
   \rang{\hat{P}}=\rang{A(I-P)}^\perp \\
   \rang{I-\hat{P}} = \rang{AP}^\perp.
   \end{gathered}
\end{equation}
\end{lemma}
\Proof For $f\in \dom(D)$ and $g\in \Lp(0,1)\otimes \C^2$,
\begin{eqnarray*}
   \dotp{ADf,g} & = & -\int_0^1 \odtp{Af''(x),g(x)} \ud x  \\
   & = & \left. \odtp{APf',g}\right|_1^0 + \int_0^1 \odtp{Af'(x),g'(x)} \ud x.
\end{eqnarray*}
We ought to find a complex $2 \times 2$ matrix $B$ and impose
boundary conditions on $g$, for
\begin{eqnarray*}
   \dotp{f,(AD)^\ast g} & = & -\int_0^1 \odtp{f(x),Bg''(x)} \ud x \\
   & = & \left. \odtp{B^\ast f,g'} \right|_1^0+ \int_0^1 \odtp{B^\ast f'(x)
   ,g'(x)} \ud x\\
   & = & \left. \odtp{B^\ast (I-P)f,g'} \right|_1^0 +
   \int_0^1 \odtp{B^\ast f'(x),g'(x)} \ud x .
\end{eqnarray*}
and
\[
\dotp{ADf,g}=\dotp{f,(AD)^\ast g}.
\]
This must be true in particular for all $f$ and $g$ with compact
support in $(0,1)$, so clearly $B=A^\ast$.

Let the boundary conditions for $(AD)^\ast$ be given by (\ref{e1})
where $\hat{P}=\hat{P}^2$ is a non-necessarily orthogonal
projection on $\C^2$, we show (\ref{e25}). If $f,g$ are smooth
functions supported in $[0,1)$, then
\[
\odtp{APf'(0),(I-\hat{P}) g(0)} = \odtp{A(I-P)f(0),\hat{P} g'(0)}
\]
where $f(0)$, $f'(0)$, $g(0)$ and $g'(0)$ are arbitrary vectors in
$\C^2$. If $f'(0)=0$, the right hand side should vanish for all
$f(0),g'(0) \in \C^2$, so that
\[
\rang{\hat{P}}=\rang{A(I-P)}^\perp.
\]
If $f(0)=0$, the left hand side should vanish for all $f'(0),g(0) \in
\C^2$, so that
\[
\rang{I-\hat{P}} = \rang{AP}^\perp.
\]
Since $A$ is non-singular these two spaces are one dimensional. \fin

\gap

\begin{corollary} \label{t19}
$AD$ is self-adjoint, if and only if $A$ is real and diagonal.
\end{corollary}
\Proof Using the notation of lemma \ref{t18}, $AD$ is
self-adjoint, if and only if
\[
  A=A^\ast \qquad \mathrm{and} \qquad P=\hat{P}.
\]
The latter occurs, if and only if
\[
   A\begin{vect} 1 \\ 0 \end{vect} \perp \begin{vect} 0 \\ 1
   \end{vect}.
\]
These conditions ensure $A$ real and diagonal. \fin

\gap

We now show that due to the boundary conditions we have chosen,
\[
  \nume{AD}=\C
\]
for a large family of non-diagonal matrices $A$. This prevents us
from employing the theory of sectorial sesquilinear forms in order
to find the spectrum.

\begin{theorem} \label{t17}
Let $A$ be a non-singular matrix.
\begin{itemize}
\item[a)] If $A$ is an upper triangular matrix (that is $c = 0$), then
\[
   \overline{\nume{AD}}=\{ rz : r\in [0,\infty), z\in \overline{\nume{A}}\}.
\]
\item[b)] If $A$ is not an upper triangular matrix (that is $c\not=0$),
then
\[
   \nume{AD}=\C.
\]
\end{itemize}
\end{theorem}
\Proof Since $0$ is always an eigenvalue of $AD$ (cf. section
\ref{s3}), then $0\in \nume{AD}$. For $f\in \dom(D)$,
\begin{eqnarray}
   \dotp{ADf,f} & = & -\int_0^1 \odtp{Af''(x),f(x)} \ud x \nonumber \\
   & = & \left. \odtp{Af',f} \right| _1^0 + \int_0^1 \odtp{Af'(x),f'(x)} \ud
   x \nonumber \\
    & = & \left.  \odtp{A\begin{vect} \phi' \\ \gam' \end{vect},
   \begin{vect} \phi \\ \gam \end{vect}} \right| _1^0
   + \int_0^1 \odtp{Af'(x),f'(x)} \ud x \nonumber \\
   & = & \left. \odtp{A\begin{vect} \phi' \\ 0 \end{vect},
   \begin{vect} 0 \\ \gam \end{vect}} \right|_1^0 +\int_0^1 \odtp{Af'(x),f'(x)}
   \ud x. \label{e13}
\end{eqnarray}

\underline{Case a)}: call
\[
   \Phi:= \{ rz : r\in [0,\infty), z\in \overline{\nume{A}}\}.
\]
Then $\Phi$ is a convex set and
\begin{eqnarray*}
   \Phi & = & \overline{\{ rz : r\in [0,\infty), z\in \nume{A}}\} \\
   & = &\overline{\{ \odtp{Av,v}: v\in \C^2 \}}.
\end{eqnarray*}
If $c=0$,
\[
   \left. \odtp{A\begin{vect} \phi' \\ 0 \end{vect},
   \begin{vect} 0 \\ \gam \end{vect}} \right|_1^0 =
   \left. \odtp{\begin{vect} a \phi' \\ 0 \end{vect},
   \begin{vect} 0 \\ \gam \end{vect}} \right|_1^0 =0
\]
so that
\[
   \dotp{ADf,f}= \int_0^1 \odtp{Af'(x),f'(x)} \ud x.
\]
This and the fact that $\Phi$ is closed and convex, yield
\[
    \overline{\nume{AD}} \subseteq \Phi.
\]

In order to prove the reverse inclusion, let $v\in \C^2$ be such
that $|v|=1$ and let
\[
   z:=\dotp{Av,v}\in \nume{A}.
\]
For all $t\geq 5$, let
\[
   \psi_t(x):= \left\{ \begin{array}{lcl} \frac{1-\cos(\pi t
   x/2)}{\sqrt{4-10/t}} & \mathrm{if} & 0\leq x \leq 2/t \\
   \frac{2}{\sqrt{4-10/t}} & \mathrm{if} & 2/t \leq x \leq 1-2/t
   \\ \frac{1-\cos(\pi t (x-1)/2)}{\sqrt{4-10/t}} & \mathrm{if} &
   1-2/t \leq x \leq 1. \end{array}\right.
\]
Then $\psi_t(0)=\psi_t(1)=\psi_t'(0)=\psi_t'(1)=0$,
\[
   \int_0^1 |\psi_t(x)|^2 \ud x =1 \qquad \mathrm{and} \qquad
   \int_0^1 |\psi_t'(x)|^2 \ud x =\frac{\pi^2t^2}{8t-20}.
\]
Let $f_t:=v \psi_t \in \dom(D)$. By construction $\|f_t\|=1$ and
\begin{eqnarray*}
   \dotp{ADf_t,f_t} & = & \int_0^1 \odtp{Af_t'(x),f_t'(x)} \ud x
   \\ & = & \odtp{Av,v} \int_0^1 |\psi_t'(x)|^2 \ud x \\
   & = & \frac{z\pi^2t^2}{8t-20}.
\end{eqnarray*}
Thus by taking $t\to\infty$, from the fact that $0\in\nume{AD}$
and since $\nume{AD}$ is convex, we gather
\[
   \overline{\nume{AD}} \supseteq \Phi.
\]

\underline{Case b)}: now $c\not=0$. Let $z$ be a fixed non-zero
complex number. Our aim is to find functions $f_\eps \in \dom(D)$
parameterized by $\eps>0$, such that $\|f_\eps\|=1$ and
$\dotp{Af_\eps,f_\eps}$ is close to $z$ for small $\eps$.

For $0<\eps<1/2$, let
\[
   \phi_\eps(x):= \left\{ \begin{array}{lcl}
   \frac{\eps}{c\pi}\sin(x\pi/\eps)
   & \mathrm{if} & 0\leq x \leq \eps/2 \\
   \frac{\eps}{2c\pi}[1-\cos(2\pi(x/\eps-1))] & \mathrm{if} & \eps/2 \leq x \leq
   \eps
   \\ 0 & \mathrm{if} &
   \eps \leq x \leq 1. \end{array}\right.
\]
Then, straightforward computations show
$\phi_\eps(0)=\phi_\eps(1)=\phi_\eps'(1)=0$,
$\phi'_\eps(0)=c^{-1}$,
\[
   \int_0^1 |\phi_\eps(x)|^2 \ud x=\frac{11\eps^3}{16c^2\pi^2}
   \qquad \mathrm{and}
   \qquad \int_0^1 |\phi_\eps'(x)|^2 \ud x=\frac{\eps}{2c^2}.
\]
For all $\eps>0$ small enough, we define the required test
function $f_\eps$ as
\[
   f_\eps(x):= \begin{vect} z \phi_\eps(x) \\
   \alp(\eps)\end{vect}
\]
where
\[
  \alp(\eps) := \sqrt{1-|z|^2\|\phi_\eps\|^2}=
  \sqrt{1-\frac{11|z|^2\eps^3}{16c^2\pi^2}}
\]
is independent of $x$. By construction $f_\eps\in \dom (D)$ and
\[
   \|f_\eps\|^2=\|z\phi_\eps\|^2+ \alp(\eps)^2 =1.
\]
According to (\ref{e13}),
\begin{eqnarray*}
   \dotp{ADf_\eps,f_\eps} & = &  \odtp{A\begin{vect} z\phi_\eps'(0)
   \\ 0 \end{vect},
   \begin{vect} 0 \\ \alp(\eps) \end{vect}}
    +\int_0^1 \odtp{Af'(x),f'(x)}
   \ud x \\
   & = & z\alp(\eps) +
   \dotp{A\begin{vect} 1 \\ 0 \end{vect},\begin{vect} 1 \\ 0
   \end{vect}} \int_0^1 |z|^2|\phi_\eps'(x)|^2 \ud x \\
   & = & z\alp(\eps) +
   \frac{a\eps}{2c^2}|z|^2.
\end{eqnarray*}

Since $\alp(\eps)\to 1$ as $\eps\to 0$, the above shows
$\dotp{ADf_\eps,f_\eps}\to z$ as $\eps\to 0$, so that $z$ is an
accumulation point of $\nume{AD}$. By moving $z\in\C$, any complex
number is accumulation point of $\nume{AD}$. Since $\nume{AD}$ is
convex, the only possibility for $\nume{AD}$ is to be the whole
complex plane. \fin

\gap


\section{The resolvent of $AD$} \label{s3}

In this section we show that the resolvent of $AD$ is compact for
all non-singular $A$. In general it is false that the product of a
bounded operator and an operator whose resolvent is compact has
compact resolvent, however if we know in addition that the
spectrum of the product is not the whole of $\C$, then the
assertion is true.

We first show that the resolvent of $D$ is compact by making use
of its self-adjointness. Since the constant function
\[
   f_0\equiv \begin{vect} 0 \\ 1 \end{vect}
\]
is in $\dom(D)$ and $ADf_0$ vanishes,
\[
0\in\spec{AD}.
\]

\begin{proposition} \label{t13}
If $A$ is a diagonal matrix, then
\[
   \spec{AD}=\{ a_-\pi^2k^2, a_+\pi^2k^2 \}_{k=0}^\infty.
\]
The zero eigenvalue is always non-degenerate and all the remaining
eigenvalues are of multiplicity no greater than 2.
\end{proposition}
\Proof Let $f_0\in
\dom(D)$ be as above. For all $n=1,2,\ldots$, let
\begin{equation} \label{e2}
f_{2n-1}(x) := \sqrt{2} \begin{vect} \sin(\pi n x) \\ 0 \end{vect}
\quad \mathrm{and} \quad f_{2n}(x) := \sqrt{2} \begin{vect}0 \\
\cos(\pi n x)
\end{vect}.
\end{equation}
Then  $f_{k}\in \dom(D)$,
\[
ADf_{2n-1}=(a_{+} \pi^2 n^2) f_{2n-1}, \qquad ADf_{2n}=(a_{-}
\pi^2 n^2) f_{2n}
\]
and $\{f_k\}_{k=0}^{\infty}$ is a complete orthonormal set in
$\Lp(0,1)\otimes \C^2$. \fin

According to corollary \ref{t19} and the above proposition,
$D=D^\ast\geq 0$ and
\[
   \spec{D}=\{\pi^2 k^2 \}_{k=0}^{\infty}.
\]
Since the eigenfunctions $\{f_k\}_{k=0}^{\infty}$ form a complete
orthonormal set, the resolvent of $D$ is compact.

\gap

Let us now rule out the possibility $\spec{AD}=\C$.

\begin{lemma} \label{t20}
For any non-singular $A\in \C$,
\[
  \spec{AD}\not=\C.
\]
\end{lemma}
\Proof Fix the matrix $A$. Since
\[
   AD-\lam = A(D-\lam A^{-1}),
\]
the complex number $\lam \in \spec{AD}$, if and only if
\[
   0\in\spec{D-\lam A^{-1}}.
\]

Let $H(\lam):=D-\lam A^{-1}$. Then the family of operators
$H(\lam)$ with domain $\dom(D)$ independent of $\lam$ is a
holomorphic family of type (A) for all $\lam \in \C$. Since 0 is a
non-degenerate isolated eigenvalue of $H(0)=D$ and $A^{-1}$ is
bounded, there exist an open neighbourhood $0\in U \subset \C$
such that $H(\lam)$ has a non-degenerate isolated eigenvalue,
(denoted by $\mu(\lam)$) close to 0 for all $\lam \in U$ and
$\mu(\lam)$ is a complex valued holomorphic function in $U$ (cf.
\cite[th.XII.8]{res4}).

If there exists some $\lam_0\in U$ satisfying $\mu(\lam_0)\not=0$,
then $0\not \in \spec{H(\lam_0)}$ so that $\lam_0\not \in
\spec{AD}$. Hence, in order to show that $\spec{AD}\not = \C$, it
is enough to show that $\mu \not \equiv 0$.  For this we find the
first coefficients in the Rayleigh-\Schrodinger series expansion
of $\mu$ about $0$. Let
\[
    \mu(\lam)=\mu_0+\mu_1 \lam + \mu_2 \lam^2 + \ldots \qquad
    \qquad \lam \in U.
\]
Since $\mu(0)=0$, $\mu_0=0$. Since $\|f_0\|=1$ and
$H(0)f_0=Df_0=0$, we compute directly $\mu_1$ (cf. \cite[remark
2.2, p.80]{kat}) by
\begin{eqnarray*}
   \mu_1 & = & \dotp{A^{-1}f_0,f_0} \\
   & = & \dotp{A^{-1}\begin{vect} 0 \\ 1 \end{vect},\begin{vect} 0 \\ 1
   \end{vect}}.
\end{eqnarray*}
If $A$ is such that $a\not = 0$,
\[
    \dotp{A^{-1}\begin{vect} 0 \\ 1 \end{vect},\begin{vect} 0 \\ 1
   \end{vect}}\not = 0
\]
so that $\mu_1$ does not vanish and hence $\mu\not \equiv 0$.

Let $A$ be such that $a=0$. Then $\mu_1=0$ so we compute $\mu_2$.
Let $f_k$ be the eigenfunctions of $D$ as in (\ref{e2}) so that
$\| f_k \|=1$ for all $k=1,2,\ldots$. Let
$\lam_{2n-1}=\lam_{2n}:=\pi^2n^2$ for all $n=1,2,\ldots$ so that
\[
   H(0)f_k=Df_k=\lam_k f_k.
\]
Then (cf. \cite[remark 2.2, p.80]{kat})
\[
   -\mu_2=\sum _{k=1}^\infty
   \frac{\dotp{A^{-1}f_0,f_k}\dotp{A^{-1}f_k,f_0}}{\lam_k}.
\]
We compute each term in the series. Since $a=0$ and $A$ is
invertible, then $b$ and $c$ do not vanish and
\[
   A^{-1}=\begin{matr2} -d/(bc) & 1/c \\ 1/b & 0 \end{matr2}.
\]
Hence
\[
   \dotp{A^{-1}f_0,f_k}= \int_0^1 \dotp{A^{-1} \begin{vect} 0 \\ 1
   \end{vect}, f_k(x)} \ud x =
   \int_0^1 \dotp{\begin{vect} 1/c \\ 0 \end{vect} ,f_k(x)}
   \ud x,
\]
so that
\[
   \dotp{A^{-1}f_0,f_{2n}}= \sqrt{2} \int_0^1 \dotp{\begin{vect} 1/c \\
   0 \end{vect},\begin{vect}0 \\ \cos(\pi n x) \end{vect}} \ud x = 0
\]
and
\begin{eqnarray*}
   \dotp{A^{-1}f_0,f_{2n-1}}& = & \sqrt{2} \int_0^1 \dotp{\begin{vect}1/c \\
   0 \end{vect},\begin{vect}\sin(\pi n x) \\ 0 \end{vect}} \ud x \\
     & = &  \sqrt{2}/c \int_0^1 \sin(\pi n x) \ud x  \\ & = &
   \left\{ \begin{array}{lcl} 0 & \mathrm{if} & n=2m \\
   2\sqrt{2}/(c\pi n) & \mathrm{if} & n=2m-1 \end{array} \right.
\end{eqnarray*}
for $m$ integer and $n=1,2,\ldots$. On the other hand
\[
   \dotp{A^{-1}f_k,f_0}= \int_0^1 \dotp{A^{-1} f_k(x),\begin{vect} 0 \\ 1
   \end{vect}} \ud x,
\]
so that
\begin{eqnarray*}
   \dotp{A^{-1}f_{2n},f_0}& = & \sqrt{2} \int_0^1 \dotp{A^{-1} \begin{vect}0 \\
   \cos(\pi n x) \end{vect},\begin{vect}0 \\ 1 \end{vect}} \ud x \\
     & = &  \sqrt{2} \int_0^1 \dotp{\begin{vect} \cos(\pi n x)/c \\ 0
    \end{vect}, \begin{vect}0 \\ 1 \end{vect}} \ud x  = 0
\end{eqnarray*}
and
\begin{eqnarray*}
   \dotp{A^{-1}f_{2n-1},f_0} & = & \sqrt{2} \int_0^1 \dotp{A^{-1}
     \begin{vect} \sin(\pi n x) \\
   0 \end{vect},\begin{vect}0 \\ 1 \end{vect}} \ud x \\
     & = &  \sqrt{2} \int_0^1 \dotp{\begin{vect} d\sin(\pi n x)/(bc)
     \\ \sin(\pi n x)/b \end{vect}, \begin{vect}0 \\ 1 \end{vect}} \ud
     x  \\ & = & \sqrt{2}/(b) \int_0^1 \sin(\pi n x) \ud x  \\ & = &
   \left\{ \begin{array}{lcl} 0 & \mathrm{if} & n=2m \\
   2\sqrt{2}/(b\pi n) & \mathrm{if} & n=2m-1 \end{array} \right.
\end{eqnarray*}
for $m$ integer and $n=1,2,\ldots$. This yields
\[
  \dotp{A^{-1}f_0,f_{k}}\dotp{A^{-1}f_k,f_{0}} =
  \left\{ \begin{array}{lcl} 0 & \mathrm{if} & k\not=4m-3 \\
   8/(bc\pi^2 n^2) & \mathrm{if} & k=4m-3 \end{array} \right.
\]
for $m=1,2,\ldots$. Thus
\[
   -\mu_2 = \sum_{m=1}^\infty \frac{8}{bc\pi^4(2m-1)^4}
    \not = 0
\]
so that $\mu \not\equiv 0$ as we required. \fin

\gap

\begin{theorem} \label{t1}
For all $z\not \in \spec{AD}$, the resolvent $(AD-z)^{-1}$ is
compact.
\end{theorem}
\Proof Since $D$ is non-negative and it has compact resolvent,
\[
AD+A=A(D+1)
\]
has a compact inverse. Let $z \not \in \spec{AD}$, then
\begin{eqnarray*}
AD-z& = & AD +A -A -z
   \\ & = & \left(I-(A+z)(AD+A)^{-1}\right)(AD+A).
\end{eqnarray*}
Hence
\[
   (AD+A)^{-1}  = (AD-z)^{-1}(I-(A+z)(AD+A)^{-1}),
\]
so that
\begin{eqnarray*}
   (AD-z)^{-1} & = & (AD-z)^{-1}(A+z)(AD+A)^{-1}+(AD+A)^{-1} \\
   & = & ((AD-z)^{-1}(A+z)+1)(AD+A)^{-1}.
\end{eqnarray*}
Thus $(AD-z)^{-1}$ is compact as needed. \fin

Theorem \ref{t1} shows that the spectrum of $AD$ consists entirely
of isolated eigenvalues of finite multiplicity. Since the
eigenvalue problem $AD f=\lam f$ is a constant coefficient system
of second order ordinary differential equations, due to the fact
that we have a combination Dirichlet and Neumann boundary
condition at both ends of the interval, the multiplicity of each
eigenvalue is never greater than 2.

\gap


\section{Asymptotics of the resolvent} \label{s4}

We now investigate the asymptotic behaviour of the resolvent norm
of $AD$.  The results we discuss in this section are connected
with the stability of the heat semigroup $\e^{-ADt}$. They are
also relevant from the computational point of view and they are
closely related to both local and global stability of the spectrum
(cf. \cite{AD1}, \cite{D2}, \cite{T1} and the reference therein).
The present approach is motivated by analogous reports on
non-self-adjoint \Schrodinger operators in \cite{B1}, \cite{D2}
and \cite{D3}.

Let
\[
   J:= \begin{matr2} 1 & 0 \\ 0 & -1 \end{matr2}.
\]
Below and elsewhere we will denote by $\tilde{D}:=JD$. According
to lemma \ref{t19}, $\tilde{D}= \tilde{D}^\ast$. According to
lemma \ref{t13},
\[
   \spec{\tilde{D}}=\{\pm \pi^2 n^2\}_{n=0}^\infty
\]
each eigenvalue being of multiplicity 1. We will employ part b) of
the following theorem in the proof of theorem \ref{t12}-b).

\begin{theorem} \label{t2}
Assume that there exists a non-singular diagonal matrix $B$ such that
$B^{-1}AB=(B^{-1}AB)^{\ast}>0$. Then
\begin{itemize}
\item[a)] $AD$ is similar to a non-negative self-adjoint operator.
\item[b)] $A\tilde{D}$ is similar to a self-adjoint operator whose
numerical range is the whole real line.
\end{itemize}
\end{theorem}
\Proof Let $C:=B^{-1}AB$ so that $C=C^\ast>0$. Since diagonal
matrices commute with the boundary conditions, $AD$ is similar to
$CD$. For the same reason and since diagonal matrices also commute
with $J$, $A\tilde{D}$ is similar to $C\tilde{D}$.

By hypothesis, the square root
$C^{1/2}=(C^{1/2})^\ast>0$. Then
\[
CD=C^{1/2}(C^{1/2}DC^{1/2})C^{-1/2}=C^{1/2}KC^{-1/2},
\]
where
\begin{gather*}
  K= C^{1/2}DC^{1/2} \\
  \dom(K)=\{ f\in \Lp(0,1)\otimes \C^2 \,:\,C^{1/2}f\in \dom(D)
  \},
\end{gather*}
so that $CD$ is similar to $K$. Since $D=D^\ast\geq 0$, then
$K=K^\ast\geq 0$.

Analogously  $C\tilde{D}$ is similar to
\[
   \tilde{K}:= C^{1/2}\tilde{D} C^{1/2}
\]
where $\dom (\tilde{K}) =\dom (K)$. Since
$\tilde{D}=\tilde{D}^\ast$, then $\tilde{K}=\tilde{K}^\ast$.
Furthermore, since
\[
   \nume{\tilde{D}}=\R
\]
and
\[
   \dotp{\tilde{K}f,f}=\dotp{\tilde{D}C^{1/2}f,C^{1/2}f},
\]
the numerical range of $\tilde{K}$ is the whole real line. \fin

Let $A$ be as in the hypothesis. The similarity to a self-adjoint
operator ensures the existence of a constant $k_A\geq 1$ such that
\[
   \|(AD-z)^{-1}\| \leq \frac{k_A}{\dist{z,[0,\infty)}} \qquad z\not \in
   \spec{AD}
\]
and
\[
   \|(A\tilde{D}-z)^{-1}\| \leq \frac{k_A}{\dist{z,\R}} \qquad z\not \in
   \spec{A\tilde{D}}.
\]
These identities show that although the numerical range of $AD$
and $A\tilde{D}$ are in general the whole complex plane, the
eigenvalues of these operators are stable in the sense of
\cite{T1}.

\gap

If we assume the weaker condition $C+C^\ast>0$, we show how to
recover part of the above estimate. We start with a preliminary
lemma.

\begin{lemma} \label{t22}
Let $A$ be such that
$
   \nume{A}\subset \{ \Re(z)>0 \}.
$
Then \linebreak $\spec{AD} \subset \{\Re(z)\geq 0\}$ and there
exists $k>0$ independent of $z$, such that
\begin{equation} \label{e26}
   \| (AD-z)^{-1} \| \leq \frac{k}{|z|} \qquad \qquad \Re(z)<0.
\end{equation}
\end{lemma}
\Proof Let $r>0$ and let $z\not \in [0,\infty)$. Then
\begin{align*}
A D-z= & A(D-zA^{-1}) \\
   = & A\left[(D-r z)+z(r -A^{-1})\right] \\
    = & A\left[ 1+(r -A^{-1} )z(D-r z)^{-1} \right](D-r z).
\end{align*}
Therefore $z\not \in \spec{A D}$, whenever
\begin{equation} \label{e4}
   \left\|(r -A^{-1})z(D-r z)^{-1} \right\| <1.
\end{equation}
We show that there is  always $r>0$  independent of $z$, such that
this holds for all $\Re(z)<0$.

Since $D\geq 0$ and $0\in \spec{AD}$,
\[
   \|(D-r z)^{-1}\| = \frac{1}{r|z|}.
\]
Thus
\[
   \left\|(r -A^{-1} )z(D-r z)^{-1} \right\| \leq
   \|1-r^{-1}A^{-1} \|.
\]
The hypothesis we imposed on $A$ is equivalent to saying
\[
   A+A^\ast>0,
\]
then
\[
   A^{-1}+(A^{-1})^\ast = A^{-1}(A^{\ast}+A)(A^{-1})^\ast>0.
\]
For all $v\in \C^2$,
\begin{eqnarray*}
   \|(I-r^{-1} A^{-1})v\|^2 & = & \dotp{\left(I-r^{-1}
    (A^{-1}+(A^{-1})^\ast) + r^{-2} (A^{-1})^\ast A^{-1}\right)v,v} \\
   & = & |v|^2 - r^{-1}\dotp{\left(A^{-1}+(A^{-1})^\ast +
   r^{-1} (A^{-1})^\ast A^{-1}\right)v,v}.
\end{eqnarray*}
Hence there exists a constant $k_0>0$ independent of $r$ (and
$z$), such that
\[
  \|I-r^{-1} A^{-1}\| < 1-r^{-1} k_0
\]
when $r$ is large enough. For such an $r$, identity (\ref{e4})
holds for any $\Re(z)<0$. This shows that $\spec{AD}$ must be
enclosed in the right hand plane. Furthermore
\begin{eqnarray*}
   \|(A D-z)^{-1} \| & \leq & \|A^{-1}\| \| (D-zA^{-1})^{-1} \| \\
   & \leq & \|A^{-1}\| \|(D-r z)^{-1}\|
   \left\| \left(1+z(D-r z)^{-1}(r -A^{-1}) \right)^{-1}\right\| \\
   & \leq & \frac{\|A^{-1}\|}{r |z| } \sum_{l=0}^{\infty}
   \|z(D-r z)^{-1}(r -A^{-1} ) \|^l \\
   & \leq & \frac{\|A^{-1}\|}{r |z| } \sum_{l=0}^{\infty}
   \|1-r^{-1}A^{-1}\|^l \\
   & \leq & \frac{k}{|z|}
\end{eqnarray*}
so (\ref{e26}) is also proven. \fin

Below and elsewhere we denote by $\Ome$ the set of non-singular
diagonal matrices and
\[
   S(\alp,\bet):= \{ z\in \C: \alp \leq \arg (z) \leq \bet \}
   \qquad \qquad \alp\leq \bet.
\]

\begin{theorem} \label{t21}
If there exists $B\in \Ome$ such that
\[
   \nume{B^{-1}A B} \subset S(\alp,\bet) \qquad \qquad
   \bet-\alp<\pi,
\]
then $\spec{AD} \subset S(\alp,\bet)$ and for any small enough
$\eps>0$ there exists $k_\eps>0$ independent of $z$,  such that
\[
   \| (AD-z)^{-1} \| \leq \frac{k_\eps}{|z|} \qquad \qquad \qquad z\not\in
   S(\alp-\eps, \bet+\eps).
\]
\end{theorem}
\Proof Let $C:=B^{-1}AB$, so that
\[
   \nume{C} \subset S(\alp,\bet).
\]
Since $B$ commutes with the boundary conditions, $AD$ is similar
to $CD$ and so it is enough to show the theorem for $CD$. Now, for
all $-(\alp+\pi/2)<\tet<\pi/2-\bet$
\[
   \nume{\e^{i\tet}C} \subset \{ \Re(z)>0 \},
\]
so we just have to apply lemma \ref{t22} to $\e^{i\tet}C$. \fin

The constant $k_\eps$ of this theorem is in general strictly
greater than $1$, therefore this is weaker than the similar
condition for m-sectorial operators in \cite[p.279]{kat}.

\gap

If $A$ is triangular, the hypothesis of the above theorem does not
necessarily hold. For instance if
\[
   A=\begin{matr2} a & 0 \\ 1 & a \end{matr2} \qquad \qquad a>0,
\]
then
\[
   \nume{A}=\{ a+z :|z|<1/2\}
\]
and so for small $a$ the numerical range contains the origin.
Nonetheless by using a similarity transformation and an
approximation argument, we can show positivity of the spectrum
whenever both of the eigenvalues of $A$ are positive ($a>0$ in our
example). The conclusion about the spectrum of the following
result will be improved in theorem \ref{t11}.

\begin{corollary} \label{t8}
Let $A$ be either upper or lower triangular. If $a\geq d
>0$, then
\begin{equation*}
   \spec{AD} \subset [0,\infty)
\end{equation*}
and for all $\eps>0$ there exists $k_\eps>0$ independent of $z$,
such that
\[
   \|(AD-z)^{-1} \| <\frac{k_\eps}{|z|}
\]
for all $z \not \in S(-\eps,\eps)$.
\end{corollary}
\Proof If $A$ is upper triangular the proof is similar so let us
assume that
\[
   A=\begin{matr2} a & 0 \\ c & d \end{matr2}.
\]
Let
\[
A(r):=\begin{matr2} 1 & 0 \\ 0 & r \end{matr2} A
   \begin{matr2} 1 & 0 \\ 0 & r^{-1} \end{matr2}=
   \begin{matr2} a &  0 \\ rc & d \end{matr2}.
\]
Then $AD$ is similar to $A(r)D$ for all $r\not=0$. Put
\[
C(r):=A(r) + A(r)^\ast = \begin{matr2} 2a & r\overline{c} \\ rc  &
2d \end{matr2}.
\]
Then $C(r)=C(r)^\ast$. The eigenvalues of $C(r)$ are
\[
    a+d \pm \sqrt{(a-d)^2 +  r^2 |c|^2},
\]
thus for small $r>0$, $C(r)>0$. The numerical range of $A(r)$ is
an ellipse with focus at $a,d$ and principal axis in the vertical
direction of the order of $r$. By taking $r\to 0$, theorem
\ref{t21} completes the proof. \fin

\gap

If $A$ is as in the hypothesis of corollary \ref{t8}, there does
not exist $B\in \Ome$ such that $B^{-1}AB=(B^{-1}AB)^\ast$ or
$B^{-1}(AJ)B=(B^{-1}(AJ)B)^\ast$ so theorem \ref{t2} is not
applicable. We show that at least in one case $AD$ fails to be
similar to self-adjoint.

\begin{theorem} \label{t24}
Let
\[
   A=\begin{matr2} a & 0 \\ 1 & a \end{matr2} \qquad \qquad a>0.
\]
Let $\eps>0$ and $z(r):= 4a\pi^2 r^2\pm i\eps$. Then there exists
a constant $k_{\eps}>0$ independent of $r$, such that
\[
   \|(AD-z(r))^{-1}\| > k_{\eps} r^{1/2}\qquad \qquad r=1,2,\ldots .
\]
\end{theorem}
\Proof Fix $\eps>0$ and let $z(r):= 4a\pi^2 r^2 - i\eps$. Without
loss of generality we can assume $r=3,4,\ldots$. Throughout the
proof the constants $l_j$ are assumed to be positive, possibly
depending upon $\eps$ but independent of $r$.  In order to show
the desired conclusion, it is enough to find $f_r \in \dom (D)$
and $l_0$, such that
\begin{equation} \label{e7}
   \frac{\| AD f_r - z(r) f_r \|}{\|f_r\|} \leq l_0  r^{-1/2}
\end{equation}
for all large enough $r$.

Let
\[
   f=\begin{vect} \phi \\ \gam \end{vect} \in\dom(D).
\]
Then
\begin{eqnarray*}
    ADf-z(r)f & = & -\begin{matr2} a & 0 \\ 1 & a \end{matr2} \begin{vect}
    \phi'' \\ \gam'' \end{vect} - z(r)\begin{vect}
    \phi \\ \gam \end{vect} \\ & = & \begin{vect} -a\phi '' -z(r)\phi \\ -\phi''
    -a \gam'' - z(r) \gam \end{vect}.
\end{eqnarray*}
Hence
\[
    \|f\|^2=\|\phi\|^2+\|\gam\|^2
\]
and
\[
   \|ADf-z(r)f\|^2  = \|a\phi''+ z(r) \phi \|^2 +\|a \gam''+z(r)
   \gam+\phi''\|^2.
\]

We now define the appropriate $f_r\in \dom(D)$ satisfying
(\ref{e7}). Let
\[
    \gam_r(x):= \cos(2\pi r x)
\]
Then $\|\gam_r\|^2=1/2$. Let
\[
    \phi_r(x):= \left\{ \begin{array}{lcl} -i\eps\cos(2\pi rx)/(4\pi^2
   r^2)
   & \mathrm{if} & x\in(1/r,1-1/r) \\ 0 & \mathrm{if} & x \not \in
   (1/(2r),1-1/(2r)) \end{array} \right.
\]
be such that $ \phi_r$ is smooth and
\begin{itemize}
\item[a)] $| \phi_r(x)|\leq \eps/(4\pi^2r^2)$  for all $x\in [0,1]$,
\item[b)] $| \phi_r'(x)|\leq l_1/r$ for all $x\not \in (1/r,1-1/r)$,
\item[c)] $| \phi_r''(x)|\leq l_2$ for all $x\not \in (1/r,1-1/r)$.
\end{itemize}
Then
\[
   \| \phi_r\|^2=\int_0^1 | \phi_r(x)|^2 \ud x \leq \eps^2/(16\pi^4r^4) \leq
   l_3r^{-4}
\]
and
\[
    \| \phi_r\|^2  \geq  \int_{1/r} ^{1-1/r} \frac{\eps^2 \cos^2(2\pi r
    x)}{16\pi^4r^4} \ud x \geq l_4r^{-4}.
\]
Hence
\[
    f_r=\begin{vect} \phi_r \\ \gam_r \end{vect} \in \dom (D)
\]
and
\begin{equation} \label{e10}
    1/2 \leq \|f_r\|^2= \| \phi_r\|^2+\| \gam_r\|^2 \leq 1
\end{equation}
for all large enough $r$. If $1/r<x<1-1/r$,
\begin{eqnarray*}
   a \phi''(x)+z(r) \phi(x) & = & a \phi(x)''+  4a\pi^2 r^2  \phi(x) - i\eps
    \phi(x)   \\
   & = & \lefteqn{ [a i \eps \cos(2\pi r x)''+ 4a\pi^2r^2 i\eps\cos(2\pi r x)
   + } \\ & & {} \qquad \qquad \qquad \qquad - \eps^2\cos(2\pi r x)]
   /(4\pi^2 r^2) \\
   & = & -\eps^2\cos(2\pi r x)/(4\pi^2r^2).
\end{eqnarray*}
Then, a) and c) yield
\begin{eqnarray}
   \|a \phi''+z(r) \phi\|^2 & = & \int _0^1 |a \phi''(x)+z(r) \phi(x)|^2
   \ud x \nonumber \\
   & \leq & \int_{1/r}^{1-1/r} l_5/r^4 \ud x +
   \int_{x\not\in[1/r,1-1/r]} l_6 +l_7/r^4 \ud x \nonumber\\
   & \leq & l_6r^{-1} + l_5r^{-4}  +l_7r^{-5} \label{e11}.
\end{eqnarray}
Also,
\begin{eqnarray*}
   a \gam''(x)+z(r) \gam(x)+ \phi''(x) & = &
   a \gam''(x)+ 4a\pi^2 r^2  \gam(x) - i\eps
    \gam(x)+ \phi''(x) \\
   & = & \lefteqn{a \cos(2\pi rx)'' + 4 a \pi^2 r^2 \cos (2\pi rx) +} \\
   & & {} \qquad \qquad \qquad \qquad -i\eps \cos (2\pi rx)
   +  \phi''(x)\\
   & = &   \phi''(x)- i\eps \cos (2\pi rx).
\end{eqnarray*}
Then for $1/r<x<1-1/r$,
\[
   a \gam''(x)+z(r) \gam(x)+ \phi''(x) = i\eps \cos (2\pi rx) -
   i\eps \cos (2\pi
   rx)=0
\]
and thus c) yields
\begin{eqnarray}
   \|a \gam''+z(r) \gam+ \phi''\|^2 & = & \int_0^1 |a \gam''(x)+z(r) \gam(x)+ \phi''(x)
   |^2 \ud x \nonumber \\
   & = & \int_{x\not\in[1/r,1-1/r]} |i\eps \cos (2\pi rx) -
    \phi''(x)|^2 \ud x \nonumber \\
   & \leq & \int_{x\not\in[1/r,1-1/r]} l_8 \ud x \nonumber\\
   & \leq & l_8r^{-1}. \label{e12}
\end{eqnarray}

In order to complete the proof for $z(r):= 4a\pi r^2 - i\eps$,
notice that (\ref{e10}), (\ref{e11}) and (\ref{e12}), show
(\ref{e7}). On the other hand, if $z(r):= 4a\pi r^2 + i\eps$ it is
enough to substitute $ \phi_r$ by $- \phi_r$ and repeat the above
computations. \fin

This result is still valid for \[ A=\begin{matr2} a & 1 \\ 0 & a
\end{matr2}.\] Indeed, it is enough to put $ \phi_r(x):=
\sin(2\pi r x)$,
\[
    \gam_r(x):= \left\{ \begin{array}{lcl} \pm i\eps\sin(2\pi rx)/(4\pi^2
   r^2)
   & \mathrm{if} & x\in(1/r,1-1/r) \\ 0 & \mathrm{if} & x \not \in
   (1/(2r),1-1/(2r)) \end{array} \right.
\]
and carry out similar calculations. Since the resolvent norm of
self-adjoint operators remains bounded in horizontal lines, the
above $AD$ can not be similar to any self-adjoint operator.

\gap

Let $\Ome_r$ be the set of all non-degenerate real diagonal
matrices. If $A$ does not satisfy the hypothesis of theorem
\ref{t21} (for instance the numerical range of $A$ is an ellipse
centered at the origin), but $A$ is ``close'' in some sense to
$\Ome_r$, an alternative to theorem \ref{t21} can be established.
We will employ this result in the proof of theorem \ref{t11}.

\begin{theorem} \label{t25}
Let there exist $B\in \Ome_r$ such that
\[
   \| AB- I\|<1.
\]
Let $\ome:=\arcsin(\|AB- I\|)$ with $0\leq \ome <\pi/2$. Then
\[\spec{AD} \subset S(-\ome,\ome)\cup S(-\pi-\ome,\ome-\pi)\] and
for any small enough $\eps>0$ there exist $k_\eps>0$ independent
of $z$, such that
\[
   \|(AD-z)^{-1} \| \leq \frac{k_\eps}{|z|}
\]
for all $z \not\in S(-\ome-\eps, \ome+\eps)\cup
S(-\pi-\ome-\eps,\ome-\pi+\eps)$.
\end{theorem}
\Proof If $\ome=0$, $A\in\Ome_r$ so the conclusion is consequence
of corollary \ref{t19}. Let $\ome>0$, let $l:= \|AB- I\|$ and put
$C:=B^{-1}\in \Ome_r$. Then $CD=(CD)^\ast$ and according to the
hypothesis $0<l<1$.

Let $z\in\C$ be such that $z \not \in S(-\ome,\ome)\cup
S(-\pi-\ome,\ome-\pi)$. Then
\begin{eqnarray*}
  (AD-z) & = & CD+(A-C)D-z \\
  & = & [I+(AB-I)CD(CD-z)^{-1}](CD-z).
\end{eqnarray*}
Since $CD$ is self-adjoint and by definition $w=\arcsin(l)$,
\begin{eqnarray}
   \| ( AB-I)CD(CD-z)^{-1} \| & \leq & l \|CD (CD-z)^{-1}\| \nonumber  \\
   & \leq & l \sup_{x\in \R} \left|\frac{x}{x-z} \right| \nonumber \\
   & \leq & \sup_{x\in \R} \frac{l}{\left| 1-\frac{z}{x}\right|}<1,
   \label{e9}
\end{eqnarray}
so that
\[
   [I+(AB-I)CD(CD-z)^{-1}]
\]
is invertible. Hence
\[
z \not \in \spec{AD}
\]
and
\begin{equation}
   (AD -z)^{-1} = (CD-z)^{-1}[I+(AB-I)CD(CD-z)^{-1}]^{-1} \label{e8}
\end{equation}
for all $\ome < |\arg(z)|\leq \pi$. This encloses $\spec{AD}$.

In order to show the second part, let \[z\not \in S(-\ome-\eps,
\ome+\eps)\cup S(-\pi-\ome-\eps,\ome-\pi+\eps),\]  for small
$\eps>0$. Then there exist a constant $l_1(\eps)>0$ independent of
$z$, such that
\[
\|(CD-z)^{-1} \| \leq \frac{l_1(\eps)}{|z|}.
\]
Also, there exist a constant $0<l_2(\eps)<1$ independent of $z$,
such that
\[
\sup_{x\in \R} \frac{l}{\left|
    1-\frac{z}{x}\right|}<l_2(\eps).
\]
These two estimates, (\ref{e9}) and (\ref{e8}) yield
\[
\|(AD-z)^{-1}\| \leq \frac{l_1(\eps) \sum_{n=0}^\infty
  l_2(\eps)^n}{|z|}=\frac{k_\eps}{|z|}. \find
\]

This shows that if $A_n\in\C^{2 \times 2}$ is a sequence of
non-singular matrices and there exists $B\in \Ome_r$ such that
\[
\|A_n - B\| \to 0
\]
as $n\to\infty$, then
\[
\spec{A_nD}\to \R.
\]

\begin{corollary} \label{t23}
Let $A$ be either upper or lower triangular. If $a,d\in \R$ and
$ad<0$, then
\begin{equation*}
   \spec{AD} \subset \R
\end{equation*}
and for all $\eps>0$ there exists $k_\eps>0$, such that
\[
   \|(AD-z)^{-1} \| <\frac{k_\eps}{|z|}
\]
for all $z \not \in S(-\eps,\eps)\cup S(-\pi-\eps,\eps-\pi)$.
\end{corollary}
\Proof It is similar to the proof of corollary \ref{t8}. Assume
without loss of generality that $b=0$ and let
\[
A(r)=\begin{matr2} 1 & 0 \\ 0 & r \end{matr2} A
   \begin{matr2} 1 & 0 \\ 0 & r^{-1} \end{matr2}=
   \begin{matr2} a &  0 \\ rc & d \end{matr2}.
\]
Then $AD$ is similar to $A(r)D$ for all $r>0$. Put
\[
C=\begin{matr2} a^{-1} & 0 \\ 0 & d^{-1} \end{matr2} \in \Ome_r,
\]
then
\[
\|A(r)C - I\| = \left\| \begin{matr2} 0 & 0 \\ rc/a  & 0
\end{matr2} \right\|= r|c/a|.
\]
Let $\ome_r:=\arcsin(r|c/a|)$. According to theorem \ref{t25}, for
all $0<r<|a/c|$
\begin{eqnarray*}
\spec{AD} & = & \spec{A(r)D} \\ & \subset & S(-\ome_r,\ome_r)\cup
S(-\pi-\ome_r,\ome_r-\pi).
\end{eqnarray*}
By taking $r$ small enough, theorem \ref{t25} yields the desired
estimate for the resolvent norm. By taking $r\to 0$, \emph{a
fortiori} $\spec{AD}\subset \R$. \fin

\gap


\section{The Hamiltonian ODE system} \label{s5}

In this section we find an entire function whose zeros coincide
with $\spec{AD}$. This is made by computing the transcendental
function associated with the $2\times 2$ system of ordinary
differential equations associated to $AD$ via standard ODE
arguments.

Let the $2\times 2$ constant coefficients second order eigenvalue
problem
\begin{gather}
   -Af''=\lam^2 f \label{e5} \\
  \begin{aligned} Pf(0)+(I-P)f'(0) & = 0   \\
  Pf(1) +(I-P)f'(1)& = 0.
  \end{aligned} \label{e15}
\end{gather}
We will say that the complex number $\lam$ is an eigenvalue of the
system (\ref{e5})-(\ref{e15}), if there exist a non-vanishing
$f\in C^\infty(0,1)\otimes \C^2$ satisfying (\ref{e5}) and the
boundary conditions (\ref{e15}).  By regularity, $\lam^2$ is an
eigenvalue of $AD$, if and only if $\lam$ is an eigenvalue of
(\ref{e5})-(\ref{e15}). Our aim is to find a holomorphic function,
denoted by $EV(x)$ below, whose zeros coincide with the
eigenvalues of (\ref{e5})-(\ref{e15}).

We proceed in the classical manner. Let the decomposition in
Jordan canonical form of $A$ be
\[
   A=:VCV^{-1}
\]
where the Jordan matrix $C$ is either
\[
   C=\begin{matr2} a_+ & 0 \\ 0 & a_-  \end{matr2} \qquad
   \mathrm{or} \qquad
   C= \begin{matr2} a_+ & 0 \\ 1 & a_+  \end{matr2}
\]
and
\[
   V:= \begin{matr2} v_1 & v_2 \\ v_3 & v_4 \end{matr2}.
\]
Then (\ref{e5})-(\ref{e15}) is equivalent to the $2\times 2$
system
\begin{gather}
  -Cg''=\lam^2 g  \label{e16} \\
  \begin{aligned}
  PVg(0)+(I-P)Vg'(0) & =  0   \\
  PVg(1) +(I-P)Vg'(1)& = 0.
  \end{aligned} \label{e17}
\end{gather}
In order to solve (\ref{e16})-(\ref{e17}), we reduce it to a first
order $4\times 4$ system as follows. For all $\lam\in\C$, let
\[
  B_\lam = \begin{matr2} 0 & I \\ -\lam^2 C^{-1} & 0 \end{matr2}
  \in \C^{4\times 4}
\]
and let
\[
   \Psi:= \left(\!\!\! \begin{array}{cccc} v_1 & v_2 & 0 & 0 \\
   0 & 0 & v_3 & v_4 \end{array}\!\!\! \right) \in \C^{2\times 4}.
\]
By regarding
\[
   \Phi= \begin{vect} g \\ g' \end{vect} \in \C^4,
\]
one sees that (\ref{e16})-(\ref{e17}) is equivalent to
\begin{gather}
  \Phi' = B_\lam \Phi \label{e18} \\
   \Psi \Phi(0) = \Psi \Phi(1) = 0. \label{e19}
\end{gather}

In order to solve (\ref{e18})-(\ref{e19}) we must find a
fundamental system of solutions. Let $e_1,e_2,e_3,e_4$ be the
standard orthonormal basis of the Euclidean space $\C^4$. A
straightforward computation show that
\[
   \exp (B_\lam x) e_j \qquad \qquad x\in [0,1], \quad j=1,2,3,4
\]
is indeed a linearly independent fundamental system for
(\ref{e18})-(\ref{e19}). Hence, $\lam$ is an eigenvalue of this
system, if and only if there exist $k_1,k_2,k_3,k_4$, such that
\begin{equation} \label{e20}
   \Phi(x) = \sum_{j=1} ^4 k_j \exp (B_\lam x) e_j
\end{equation}
is non-vanishing and satisfies the boundary conditions.

We now proceed to compute $EV(x)$. The exponential of $B_\lam x$
is given by
\[
   \exp(B_\lam x)  =  \begin{matr2} \cos(\lam C^{-1/2}x) &
   \lam^{-1}C^{1/2} \sin (\lam C^{-1/2}x) \vspace{.2in}\\
   -\lam C^{-1/2} \sin (\lam C^{-1/2} x) & \cos (\lam C^{-1/2}
   x) \end{matr2}
\]
for $x\in [0,1]$. In theorems \ref{t3} and \ref{t26} below, we
split our computation into two cases depending upon the Jordan
matrix $C$.

\begin{theorem} \label{t3}
When
\[
    C=\begin{matr2} a_+ & 0 \\ 0 & a_-  \end{matr2},
\]
$\lam$ is an eigenvalue of the system (\ref{e18})-(\ref{e19}), if
and only if $EV(\lam)=0$ for
\begin{eqnarray*}
   EV(x):=\lefteqn{\left(2\prod_{j=1}^4 v_j\right)
    \left[1-\cos \left(\frac {x}{\sqrt{a_+}}
   \right)\,
   \cos \left(\frac {x}{\sqrt{a_- }}\right) \right] + {} }  \\
   & & {} - \left[ v_1^2v_4^2\frac{\sqrt{a_+ }}{\sqrt{a_-}}+
   v_2^2 v_3^2 \frac{\sqrt{a_-}}{\sqrt{a_+ }}
   \right] \sin \left(\frac {x}{\sqrt{a_+}}\right)\,
   \sin \left(\frac {x}{\sqrt{a_- }}\right).
\end{eqnarray*}
\end{theorem}
\Proof  Notice that $EV(0)=0$. Assume $\lam\not=0$. According to
the hypothesis,
\[
    C^{1/2}=\begin{matr2} a_+^{1/2} & 0 \\ 0 & a_-^{1/2}  \end{matr2}
  \qquad \mathrm{and} \qquad
  C^{-1/2}=\begin{matr2} a_+^{-1/2} & 0 \\ 0 & a_-^{-1/2}
  \end{matr2}.
\]
Then
\[
   \exp(B_\lam x)= \left(\!\!\!
   \begin{array}{cccc} \cos \frac{\lam x}{\sqrt{a_+}} & 0 &
   \frac{\sqrt{a_+}}{\lam} \sin \frac{\lam x}{\sqrt{a_+}} & 0 \\
   0 & \cos \frac{\lam x}{\sqrt{a_-}} & 0 & \frac{\sqrt{a_-}}{\lam}
   \sin \frac{\lam x}{\sqrt{a_-}} \\ -\frac{\lam}{\sqrt{a_+}} \sin
   \frac{\lam x}{\sqrt{a_+}} & 0 & \cos \frac{\lam x}{\sqrt{a_+}}
   & 0 \\ 0 & -\frac{\lam x}{\sqrt{a_-}} \sin \frac{\lam x}{\sqrt{a_-}}
   & 0 & \cos \frac{\lam x}{\sqrt{a_-}} \end{array} \!\!\!
   \right).
\]

Let $\Phi(x)$ be a particular solution given as in (\ref{e20}),
where the complex parameters $k_j$ are to be determined. Then
\[
   \Psi \Phi(0)= \begin{vect} k_1 v_1 + k_2 v_2 \\
   k_3 v_3 + k_4 v_4\end{vect}
\]
and
\[
   \Psi \Phi(1) = \left( \!\!\! \begin{array}{l}
    k_1v_1 \cos \frac{\lam}{\sqrt{a_+}}
   + k_2 v_2 \cos \frac{\lam}{\sqrt{a_-}} + \\  \hspace{.9in}
    +k_3 v_1
   \frac{\sqrt{a_+}}{\lam} \sin \frac{\lam}{\sqrt{a_+}} +
   k_4 v_2 \frac{\sqrt{a_-}}{\lam} \sin \frac{\lam}{\sqrt{a_-}} \\
   -k_1 v_3 \frac{\lam}{\sqrt{a_+}} \sin \frac{\lam}{\sqrt{a_+}} -
   k_2 v_4 \frac{\lam}{\sqrt{a_-}} \sin \frac{\lam}{\sqrt{a_-}} +
   \\ \hspace{1.4in}
   + k_3 v_3 \cos \frac{\lam}{\sqrt{a_+}} + k_4 v_4 \cos
   \frac{\lam}{\sqrt{a_-}} \end{array} \!\!\!\right).
\]
The solution $\Phi$ satisfies the boundary conditions (\ref{e19}),
if and only if
\[
\left\{ \begin{array}{l}
  k_1 v_1 + k_2 v_2 =0 \\
  k_3 v_3 + k_4 v_4 =0 \\
   k_1v_1 \cos \frac{\lam}{\sqrt{a_+}}
   + k_2 v_2 \cos \frac{\lam}{\sqrt{a_-}}
    +k_3
   \frac{v_1\sqrt{a_+}}{\lam} \sin \frac{\lam}{\sqrt{a_+}} +
   k_4  \frac{v_2\sqrt{a_-}}{\lam} \sin \frac{\lam}{\sqrt{a_-}} =0 \\
   -k_1  \frac{v_3\lam}{\sqrt{a_+}} \sin \frac{\lam}{\sqrt{a_+}} -
   k_2 \frac{v_4 \lam}{\sqrt{a_-}} \sin \frac{\lam}{\sqrt{a_-}}

   + k_3 v_3 \cos \frac{\lam}{\sqrt{a_+}} + k_4 v_4 \cos
   \frac{\lam}{\sqrt{a_-}} =0. \end{array} \right.
\]
The determinant of this $4\times 4$ system of linear equations in
$k_j$ is precisely $EV(\lam)$. \fin

\gap

\begin{theorem} \label{t26}
When
\[
    C=\begin{matr2} a_+ & 0 \\ 1 & a_+  \end{matr2},
\]
$\lam$ is an eigenvalue of the system (\ref{e18})-(\ref{e19}), if
and only if $EV(\lam)=0$ for
\begin{eqnarray*}
   EV(x):= \left( \frac{v_2^2 v_4^2}{4a_+^3}\right) x^2 -
   \left(\det V + \frac{v_2v_4}{2a_+} \right)^2 \sin ^2\frac{x}
   {\sqrt{a_+}}.
\end{eqnarray*}
\end{theorem}
\Proof Notice that $EV(0)=0$. Assume $\lam\not=0$. One can verify
directly that
\[
    C^{1/2}=\begin{matr2} a_+^{1/2} & 0 \\
    \frac{1}{2\sqrt{a_+}} & a_+^{1/2}  \end{matr2}
  \qquad \mathrm{and} \qquad
  C^{-1/2}=\begin{matr2} a_+^{-1/2} & 0 \\
   -\frac{1}{\left(2a_+^{3/2}\right)} & a_+^{-1/2}
  \end{matr2}.
\]
Then the four $2\times2$ blocks of the matrix $\exp(B_\lam x)$ are
\[
   \cos(\lam C^{-1/2} x)  =
   \begin{matr2} \cos \frac{\lam x}{\sqrt{a_+}} & 0 \\
   \frac{\lam x}{\left(2a_+^{3/2}\right)} \sin \frac{\lam x}{\sqrt{a_+}} &
   \cos \frac{\lam x}{\sqrt{a_+}} \end{matr2},
\]
$\lam^{-1}C^{1/2} \sin(\lam C^{-1/2} x)$ equal to
\[
   \begin{matr2} \frac{a_+^{1/2}}{\lam} \sin \frac{\lam x}{\sqrt{a_+}}
   & 0 \\
   \left[ \frac{1}{2\lam \sqrt{a_+}} \sin \frac{\lam x}{\sqrt{a_+}}
   - \frac{x}{2a_+} \cos \frac{\lam x}{\sqrt{a_+}}\right]  &
   \frac{a_+^{1/2}}{\lam} \sin \frac{\lam x}{\sqrt{a_+}} \end{matr2}
\]
and $-\lam C^{-1/2} \sin(\lam C^{-1/2} x)$ equal to
\[
   \begin{matr2} -\frac{\lam }{a_+^{1/2}} \sin \frac{\lam x}{\sqrt{a_+}}
   & 0 \\
   \left[\frac{\lam}{\left(2a_+^{3/2}\right)} \sin \frac{\lam x}{\sqrt{a_+}}
   + \frac{\lam^2 x}{2a_+^2} \cos \frac{\lam x}{\sqrt{a_+}} \right]  &
   -\frac{\lam }{a_+^{1/2}} \sin \frac{\lam x}{\sqrt{a_+}} \end{matr2}.
\]

Let $\Phi(x)$ be a particular solution given as in (\ref{e20}),
where the complex parameters $k_j$ are to be determined. Then
\[
   \Psi \Phi(0)= \begin{vect} k_1 v_1 + k_2 v_2 \\
   k_3 v_3 + k_4 v_4\end{vect}
\]
and
\[
   \Psi \Phi(1)= \begin{vect} \Psi \Phi(1)_1 \\ \Psi \Phi(1)_2
   \end{vect}
\]
where
\[
   \begin{array}{l}
   \Psi \Phi(1)_1 =
    k_1\left( v_1 \cos \frac{\lam}{\sqrt{a_+}}+\frac{v_2
    \lam}{2a_+^{3/2}}\sin \frac{\lam}{\sqrt{a_+}} \right)
   + k_2 v_2 \cos \frac{\lam}{\sqrt{a_+}} + \\ \hspace{.3in}
    +k_3 \left(
   \frac{v_1\sqrt{a_+}}{\lam} \sin \frac{\lam}{\sqrt{a_+}}+
   \frac{v_2}{2\lam \sqrt{a_+}} \sin \frac{\lam}{\sqrt{a_+}}
   -\frac{v_2}{2a_+} \cos\frac{\lam}{\sqrt{a_+}} \right) +
   k_4 \frac{v_2 \sqrt{a_+}}{\lam} \sin \frac{\lam}{\sqrt{a_+}}
\end{array}
\]
and
\[
   \begin{array}{l}
   \Psi \Phi(1)_2 =
   k_1 \left(- \frac{v_3\lam}{\sqrt{a_+}} \sin \frac{\lam}
   {\sqrt{a_+}} +
   \frac{v_4\lam}{2a_+^{3/2}} \sin\frac{\lam}{\sqrt{a_+}}
    + \frac{v_4 \lam^2}{2a_+} \cos \frac{\lam}{\sqrt{a_+}}
   \right) + \\ \hspace{.3in} -
   k_2  \frac{v_4\lam}{\sqrt{a_+}} \sin \frac{\lam}{\sqrt{a_+}}
      + k_3 \left(v_3 \cos \frac{\lam}{\sqrt{a_+}} + \frac{v_4
   \lam}{2a_+^{3/2}}\sin \frac{\lam}{\sqrt{a_+}} \right) + k_4 v_4 \cos
   \frac{\lam}{\sqrt{a_+}}. \end{array}
\]
The solution $\Phi$ satisfies the boundary conditions (\ref{e19}),
if and only if
\[
\left\{ \begin{array}{l}
  k_1 v_1 + k_2 v_2 =0 \\
  k_3 v_3 + k_4 v_4 =0 \\
   \Psi \Phi(1)_1=0 \\
   \Psi \Phi(1)_2=0 \end{array} \right.
\]
A rather long but straightforward computation shows that the
determinant of this $4\times 4$ system of linear equations in
$k_j$ is $EV(\lam)$. \fin

\gap

We show that $AD$ can have non-real eigenvalues even when the
spectrum of $A$ is positive.

\begin{example} \label{t4}
Put
\[
A:= \begin{matr2} 2/5+3i/10 & 3/5-3i/10 \\ 3/20+3i/10 & 17/20-3i/10
   \end{matr2}.
\]
Then the eigenvalues of $A$ are $a_+=1$, $a_- =1/4$, and the
eigenvectors
\[
v_+ =\begin{vect} 1 \\ 1 \end{vect}, \qquad v_- =\begin{vect} 2i
\\ 1
   \end{vect}.
\]
Thus
\[
EV(x)=4i(1-\cos(x)\cos(2x))=4i(1-2\cos^3(x)+\cos(x))
\]
so that $EV(\lam)$=0, if and only if
\[
\cos(\lam) = 1 \qquad \mathrm{or} \qquad \cos(\lam)=-1/2 \pm i/2.
\]
Hence
\[
\spec{AD} = \{ 4k^2\pi^2,\, (\lam_\pm + 2k\pi)^2 \}_{k\in \Z}
\]
where $\lam_\pm=\arccos(-1/2 \pm i/2) \approx 2.02 \pm 0.53 i$.
\end{example}

\gap


\section{Real matrices} \label{s6}

In this section we explore some connections between the entries of
the matrix $A$ and the global behaviour of $\spec{AD}$ when $A
\in\R^{2\times 2}$. Alongside we discuss conditions to ensure
similarity to a self-adjoint operator. For completeness of the
picture, below and elsewhere we allow $\det(A)=0$.

Our first task is to reduce to two parameters the four that are
initially given as entries of $A$. This leads us to five different
types of matrices to deal with. For $a,d\in\R$, let
\[
\begin{array}{ccc}
   A_0:=\begin{matr2} a & 0 \\ 0 & d \end{matr2}, &
   A_1:=\begin{matr2} a & 1 \\ 1 & d \end{matr2}, &
   A_2:=\begin{matr2} a & 0 \\ 1 & d \end{matr2}, \\
   A_3:=\begin{matr2} a & 1 \\ 0 & d \end{matr2} & \mathrm{and} &
   A_4:=\begin{matr2} a & -1 \\ 1 & d \end{matr2}.\\
\end{array}
\]
We show that the $A_jD$ generate any $AD, A\in \R^{2\times 2}$ via
similarity transformations.

\begin{lemma} \label{t9}
If $A\in \R^{2\times 2}$, then $AD$ is similar to $\alp A_jD$ for
some $\alp,a,d\in \R$ and $j=0,\ldots,4$.
\end{lemma}
\Proof Let
\[
   A= \begin{matr2} \tilde{a} & b \\ c & \tilde{d} \end{matr2}.
\]
If $bc=0$, the proof is trivial. Let
\[
A(r):=\begin{matr2} 1 & 0 \\ 0 & r \end{matr2} A
   \begin{matr2} 1 & 0 \\ 0 & r^{-1} \end{matr2}=
   \begin{matr2} \tilde{a} & r^{-1}b \\ rc & \tilde{d} \end{matr2}
\]
Then, $AD$ is similar to $A(r)D$ for all $r\not=0$. If $b/c>0$,
\[
   A(\sqrt{b/c})=\begin{matr2} \tilde{a} & \sqrt{bc} \\ \sqrt{bc} & \tilde{d}
   \end{matr2} =\alp  A_1
\]
for $\alp=\sqrt{bc}$, $a=\tilde{a}/\sqrt{bc}$ and
$d=\tilde{d}/\sqrt{bc}$. If $b/c<0$,
\[
   A(\sqrt{-b/c})=\begin{matr2} \tilde{a} & \mp\sqrt{-bc} \\ \pm\sqrt{-bc}
   & \tilde{d}\end{matr2} = \pm\alp A_{4}
\]
for $\alp=\sqrt{-bc}$, $a=\pm\tilde{a}/\sqrt{-bc}$ and
$d=\pm\tilde{d}/\sqrt{-bc}$. \fin

The case $j=0$ was already described in corollary \ref{t19}.
Indeed if $ad\not=0$ then $A_0D$ is similar to a self-adjoint
operator and
\[
   \spec{A_0D}=\{an^2\pi^2,dn^2\pi^2\}_{n=0}^\infty \subset \R.
\]

\gap

\subsection{Matrix $A_1$} \ \label{ss6.1} \newline
Since $a$ and $d$ are real, $A_1=A_1^\ast$. Let $b_\pm$ be the
eigenvalues of $A_1$. Then
\[
   b_\pm=\frac{a+d \pm \sqrt{(a-d)^2+4}}{2},
\]
so that
\begin{itemize}
\item[i)] $b_+\geq b_->0$, if and only if $ad>1$ and $a,d>0$,
\item[ii)] $b_-\leq b_+<0$, if and only if $ad>1$ and $a,d<0$,
\item[iii)] $b_+$ and $b_-$ have opposite signs, if and only if
$ad<1$.
\end{itemize}

\begin{theorem} \label{t10}
The following statements are true.
\begin{itemize}
\item[a)] If $ad=1$ then $\spec{A_1D}=\C$.
\item[b)] If $ad>1$ and $a,d>0$ then $A_1D$ is
similar to a non-negative
operator so that $\spec{A_1D}\subset [0,\infty)$.
\item[c)] If $ad>1$ and $a,d<0$ then $-A_1D$
is similar to a non-negative
self-adjoint operator so that $\spec{A_1D}\subset (-\infty,0]$.
\item[d)] If $ad<1$ then $\spec{A_1D} \subset \R$.
\end{itemize}
\end{theorem}
\Proof If $ad=1$,  the matrix $A_1$ is singular so according to
lemma \ref{t27}, $A_1D$ is not a closed operator. This shows a).
Statement b) is consequence of i) and theorem \ref{t2}, and
statement c) is consequence of ii) and theorem \ref{t2}.

Let us show d). For $\eps\in \R$, let
\[
   B(\eps):= A_1+i\eps.
\]
Then
\[
   \nume{B(\eps)}\subset \{\Im(z)>0\} \qquad \qquad \eps>0
\]
and
\[
   \nume{B(\eps)}\subset \{\Im(z)<0\} \qquad \qquad \eps<0.
\]
According to theorem \ref{t21},
\[
   \spec{B(\eps)D}\subset \{\Im(z)\geq 0\} \qquad \qquad \eps>0
\]
and
\[
   \spec{B(\eps)D}\subset \{\Im(z)\leq0\} \qquad \qquad \eps<0.
\]
Since $B(\eps)D$ is a holomorphic family of type (A) in a
neighbourhood of $\eps=0$ and $B(0)=A_1$, \emph{a fortiori}
\[
   \spec{A_1D}\subset  \R. \find
\]

Although $A_1=A_1^\ast$, it is unclear to us whether $A_1D$ is
similar to self-adjoint in the latter case.

\gap

\subsection{Matrices $A_2$ and $A_3$} \ \label{ss6.2} \newline
Since the results for the matrix $A_3$ are analogous and shown in
a similar manner as for $A_2$, we will only consider the latter.

\begin{theorem} \label{t11}
The following statements are true.
\begin{itemize}
\item[a)] If $ad=0$ then $\spec{A_2D} =\C$.
\item[b)] If $ad\not=0$ then $\spec{A_2D} = \{a\pi^2n^2,d\pi^2n^2\}_{n=0}^\infty$.
\item[c)] If $ad>0$, for all $\eps>0$ there exists $k_\eps>0$ independent of $z$,
such that
\[
   \|(A_2D-z)^{-1}\| \leq \frac{k_\eps}{|z|} \qquad \qquad z\not \in
   \pm S(-\eps,\eps),
\]
where the symbol $\pm$ is chosen according to the symbol of $a$.
\item[d)] If $ad<0$,
then for all $\eps>0$ there exists $k_\eps>0$ independent of $z$,
such that
\[
   \|(A_2D-z)^{-1}\| \leq \frac{k_\eps}{|z|}
\]
for all $ z\not \in S(-\eps,\eps)\cup S(-\pi-\eps,-\pi+\eps).$
\item[e)] If $a=d\not=0$, let $\eps>0$ and $z_r=4a\pi^2r^2\pm i\eps$.
Then there exists a constant $k_\eps>0$ independent of $r$, such
that
\[
   \|(A_2D-z_r)^{-1}\|\geq k_\eps |z_r|^{1/4}
\]
for all $r=1,2,\ldots$.
\end{itemize}
\end{theorem}
\Proof If $ad=0$,  the matrix $A_2$ is singular so according to
lemma \ref{t27}, $A_2D$ is not a closed operator. This shows a).

Let us show b). If $a\not=d$, the matrix $A_2$ is diagonalizable
and
\[
   A_2=\begin{matr2} a-d & 0 \\ 1 & 1 \end{matr2}
   \begin{matr2} a & 0 \\ 0 & d \end{matr2}
   \begin{matr2} (a-d)^{-1} & 0 \\ -(a-d)^{-1} & 1
   \end{matr2}.
\]
Then, according to theorem \ref{t3},
\[
 EV(x)  =  k_0 \sin \left(\frac {x}{\sqrt{a}}\right)\,
   \sin \left(\frac {x}{\sqrt{d }}\right)
\]
where $k_0$ is constant in $x$. If $a=d$, $A_2$ is already in
Jordan form and according to theorem \ref{t26},
\[
   EV(x)= -\left[\sin\left(\frac{x}{\sqrt{a}}\right)\right] ^2.
\]
Hence in both cases
\[
    \spec{A_2D} = \{a\pi^2n^2,d\pi^2n^2\}_{n=0}^\infty.
\]

 Statements c) is consequence of corollary \ref{t8} and
statement d) is consequence of corollary \ref{t23}. For statement
e) use theorem \ref{t24} and the fact that $|z_r|$ is of order
$r^2$. \fin

\gap

\subsection{Matrix $A_4$} \ \label{ss6.3} \newline
Formally speaking, so far the spectrum of $A_jD$ for
$j=0,\ldots,3$ reproduces the spectrum of $A_j$ in the following
sense: if $A_j$ is non-degenerated and both eigenvalues of $A_j$
are positive (negative) then $\spec{A_jD}$ is non-negative
(non-positive), and if the eigenvalues are of opposite sign then
$A_jD$ possess both positive and negative spectrum. There is no
reason to expect the same for $j=4$, in fact this case is less
simple due to the way the entries of $A_4$ interact with the
boundary conditions.

The eigenvalues of $A_4$ are given by
\begin{equation} \label{e21}
   b_\pm:=\frac{a+d\pm \sqrt{(a-d)^2-4}}{2}.
\end{equation}
Then
\begin{itemize}
\item[i)] $b_+= b_-$, if and only if $|a-d|=2$. In this case $A_4$
is not a diagonalizable matrix.
\item[ii)] $b_\pm$ are real and have opposite signs, if and only if
$ad<-1$.
\item[iii)] $b_+> b_->0$, if and only if $ad>-1$,
$|a-d|>2$ and $a+d>0$.
\item[iv)] $b_-< b_+ <0$, if and only if $ad>-1$, $|a-d|>2$ and
$a+d<0$.
\item[v)] $b_\pm$ are non-real with $b_+=\overline{b_-}$, if and
only if $|a-d|<2$.
\item[vi)] $A_4$ is singular, if and only if $ad=-1$.
\end{itemize}
Motivated by this and for simplicity, we can divide the plane into
$6$ disjoint regions $R_k$,
\begin{eqnarray*}
   R_1 & := & \{(a,d)\in \R^2 : |a-d|=2 , a\not = \pm 1\},  \\
   R_2 & := & \{(a,d)\in \R^2 : ad<-1 \},  \\
   R_3 & := & \{(a,d)\in \R^2 : ad>-1, |a-d|>2, a+d>0 \},  \\
   R_4 & := & \{(a,d)\in \R^2 : ad>-1, |a-d|>2, a+d<0 \},  \\
   R_5 & := & \{(a,d)\in \R^2 : |a-d|<2\},  \\
   R_6 & := & \{(a,d)\in \R^2 : ad=-1\}.
\end{eqnarray*}
Clearly $\R^2= \bigcup R_k$. Below we establish the spectral
results for $A_4D$ separately in each region $R_k$.

\begin{figure}[t]
\begin{picture}(300,300)(55,120) \includegraphics{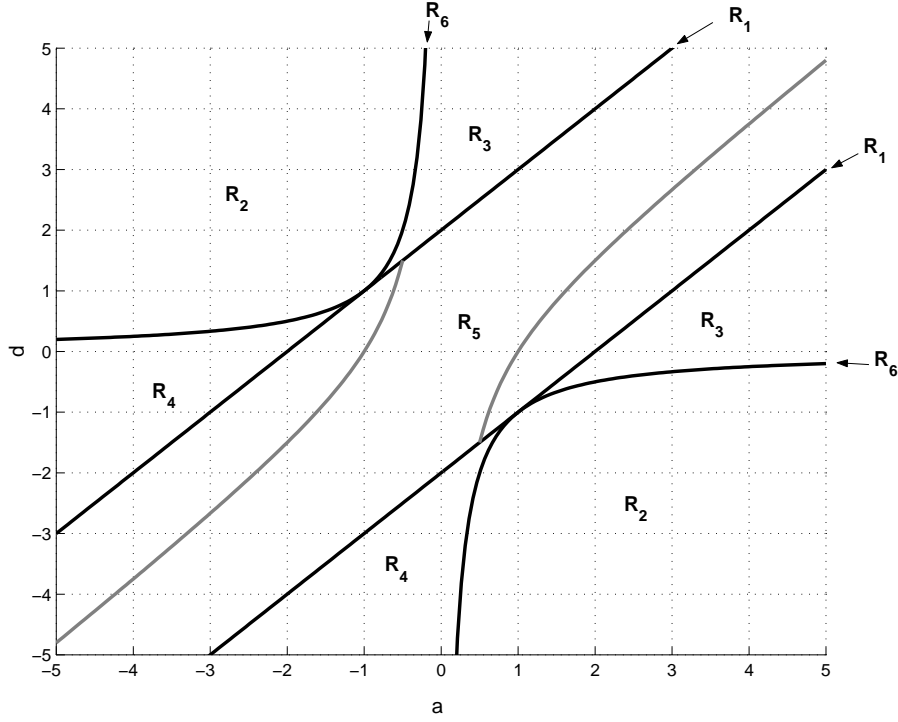}
\end{picture}
\caption{Different regions of the plane in which $\spec{A_4D}$
exhibits a similar behaviour. The grey line is $\{a^2-ad-1=0\}
\cap R_5$. See theorems \ref{t12}-\ref{t32}.} \label{f1}
\end{figure}

Two cases are similar to what we have found so far.

\begin{theorem} \label{t12} The following statements are true.
\begin{itemize}
\item[a)] If $(a,d)\in R_6$, then $\spec{A_4D}=\C$.
\item[b)] If $(a,d)\in R_2$, then $\spec{A_4D}\subset \R$ and $A_4D$
is similar to a self-adjoint operator whose numerical range is the
whole real line.
\end{itemize}
\end{theorem}
\Proof If $ad=-1$,  the matrix $A_4$ is singular so according to
lemma \ref{t27}, $A_4D$ is not a closed operator. This shows a).

Let us show b). Let $J$ be as in section \ref{s4}. Then
\[
   A_4D=(A_4J)(JD)=\begin{matr2} a & 1 \\ 1 & -d \end{matr2}
   \tilde{D}=\tilde{A}\tilde{D}.
\]
Here $\tilde{A}=\tilde{A}^\ast$ and the eigenvalues of $\tilde{A}$
are
\[
   \tilde{b}_\pm=\frac{a-d\pm \sqrt{(a+d)^2+4}}{2}.
\]
Since $ad<-1$, $\tilde{b}_\pm$ are either both positive or both
negative. If they are both positive, $\tilde{A}>0$ so that theorem
\ref{t2}-b) provides the desired conclusion. If they are both
negative apply the above argument to $-A_4D$. \fin

\gap

In order to find $\spec{A_4D}$ in $R_k$ for $k=1,3,4,5$, we ought
to rely on properties of the transcendental function $EV(x)$.
Nonetheless, theorem \ref{t30} provides some indication of what we
should expect, it bases on the observation that if both $a$ and
$d$ are positive,
\[
   A_4+A_4^\ast=\begin{matr2} 2a & 0 \\ 0 & 2d \end{matr2} >0,
\]
so by virtue of lemma \ref{t22}, $\spec{A_4D} \subset \{
\Re(z)\geq 0 \}$.

\begin{theorem} \label{t30}
If both $a$ and $d$ are positive, then
\[
   \spec{A_4D} \subset S(-\ome,\ome)
\]
where $\sin \ome =1/\sqrt{ad+1}$ for $0<\ome<\pi/2$.
\end{theorem}
\Proof The numerical range of $A_4$ is an ellipse whose foci are
$b_\pm$ and largest diameter is of length $|a-d|$. It is easy to
see that $S(-\ome,\ome)$ is the minimal sector that contains such
an ellipse. Use theorem \ref{t21} to complete the proof. \fin

Since
\[
   \begin{matr2} -1 & 0 \\ 0 & 1 \end{matr2}
   \begin{matr2} -a & -1 \\ 1 & -d \end{matr2}
   \begin{matr2} -1 & 0 \\ 0 & 1 \end{matr2}=
   -\begin{matr2} a & -1 \\ 1 & d \end{matr2}
\]
and because of diagonal matrices commute with the boundary
conditions, $\spec{A_4D} \subset -S(-\ome,\ome)$ where both $a$
and $d$ are negative. This also shows that the spectral results
for $A_4D$ are symmetric with respect to the  transformation
$(a,d)\mapsto (-a,-d)$. Below we will employ this symmetry often
without mention.

 \gap

In order to describe $\spec{A_4D}$ in $R_5$, we will make use of
the following technical result.

\begin{lemma} \label{t16}
Let $\alp\in\C$ be such that $\Re(\alp^2)\geq 0$, let $-1\leq c
\leq 1$ and let
\[
   F(x):=1-\cos(\alp x) \cos(\overline{\alp}x) - c \sin(\alp x)
   \sin(\overline{\alp}x) \qquad \qquad x\in\C.
\]
Then $F(x)$ has an infinite number of zeros in the complex plane
and
\begin{itemize}
\item[a)] if $c=-1$, then $F(x)=0$, if and only if
$\sin\left(\Re(\alp\right) x)=0$,
\item[b)] if $c=1$, then $F(x)=0$, if and only if
$\sinh\left(\Im(\alp)x\right)=0$,
\item[c)] if $-1<c<1$, then $F(x)$ only has a finite number of zeros
lying on the real and imaginary axis.
\end{itemize}
\end{lemma}
\Proof Let  $\alp=:\rho +i \mu$ so that $\rho\geq \mu>0$ and let
$x=:x_1+ix_2$ for $x_1,x_2\in\R$.

In order to show a), assume $c=-1$. Then
\begin{eqnarray*}
   |F(x)|^2 & = & \left| 1- \cos[\alp (x_1+ix_2)] \cos[\overline{\alp}
    (x_1+ix_2)] + \right.\\
   & & \qquad \qquad \qquad + \left.\sin[\alp (x_1+ix_2)]
    \sin[\overline{\alp} (x_1+ix_2)]\right| ^2 \\
   & = & 4 [\cos^2(\rho x_1) - \cosh^2(\rho x_2)]^2.
\end{eqnarray*}
Hence
\[
   F(x)=0,
\]
if and only if $\cosh(\rho x_2)=1$ and $\cos(\rho x_1)=\pm 1$.
This gives a).

Similarly for b), assume $c=1$. Then
\begin{eqnarray*}
   |F(x)|^2 & = & \left| 1- \cos[\alp (x_1+ix_2)] \cos[\overline{\alp}
    (x_1+ix_2)] + \right.\\
   & & \qquad \qquad \qquad - \left.\sin[\alp (x_1+ix_2)]
    \sin[\overline{\alp} (x_1+ix_2)]\right| ^2 \\
   & = & 4 [\cosh^2(\mu x_1) - \cos^2(\mu x_2)]^2.
\end{eqnarray*}
Hence
\[
   F(x)=0,
\]
if and only if $\cosh(\mu x_1)=1$ and $\cos(\mu x_2)=\pm 1$.

Let us show assertion c). If $x\in\R$, then
\begin{eqnarray*}
   F(x) & = & 1-\cos(\alp x)\overline{\cos(\alp x)}
   -c \sin(\alp x) \overline{\sin (\alp x)} \\
   & = & 1-|\cos(\alp x)|^2 -c |\sin(\alp x)|^2 \\
   & = & 1-\cos^2(\rho x) - c \sin^2(\rho x)  - (1+
   c) \sinh ^2(\mu x)
\end{eqnarray*}
and
\begin{eqnarray*}
   F(ix) & = & 1-\cos(-\overline{i\overline{\alp}} x)
   \cos(i \overline{\alp} x)
   -c \sin(-\overline{i\overline{\alp}} x) \sin (i \overline{\alp} x) \\
   & = &1-\overline{\cos(i\overline{\alp} x)}
   \cos(i \overline{\alp} x)
   +c \overline{\sin(i \overline{\alp} x)} \sin (i \overline{\alp} x) \\
   & = & 1-|\cos(i \overline{\alp} x)|^2 +c |
   \sin(i \overline{\alp} x)|^2 \\
   & = &1-\cos^2(\mu x) + c \sin^2(\mu x)  - (1-
   c) \sinh^2 (\rho x).
\end{eqnarray*}
Hence, if $-1<c<1$,
\[
    \lim_{x\to \pm \infty} F(x)=-\infty \qquad \mathrm{and} \qquad
    \lim_{x\to\pm \infty} F(ix)=-\infty.
\]
Since $F(x)$ is a smooth function, c) follows.

Finally let us show that $F(x)$ has a infinite number of zeros.
Suppose that $F$ only has a finite number of zeros $0,z_1,\ldots,
z_n$ where the $z_j$ repeats as many times as its order. Then
\[
   G(x)=\frac{F(x)}{x^2 \prod_{j=1}^n (x-z_j)}
\]
is an entire function with no zeros. By virtue of the Weierstrass
factorization theorem, there is an entire function $g(x)$ such
that $G(x)= \e^{g(x)}$. Then
\[
   F(x)= \left[ x^2\prod_{j=1}^n (x-z_j) \right] \e^{g(x)}
   =:p(x) \e^{g(x)} .
\]
Since it is a combination of sines and cosines, the order (cf.
\cite[p.285]{con}) in the sense of entire functions of $F(x)$ is
$\lam=1$. Thus by virtue of Hadamard's factorization theorem,
$g(x)$ is a polynomial of degree 1 in $x$ and so
\[
   F(x) = p(x) \e^{k x+ l}
\]
for suitable $k,l \in\C$. Since $p(x)$ is a polynomial, this is
clearly a contradiction, so $F(x)$ must have an infinite number of
zeros. \fin

\begin{theorem} \label{t31} Let $(a,d)\in R_5$.
\begin{itemize}
\item[a)] If $(a,d) \in \{ a^2-ad-1=0 \} \cap \{-2< a-d<0\}$, then
\[\spec{A_4D} =\left\{ -k^2\pi^2/ [\Im(b_+^{-1/2})]^2
\right\}_{k\in \Z} \subset (-\infty,0].\]
\item[b)] If $(a,d) \in \{ a^2-ad-1=0 \} \cap \{0< a-d<2\}$, then
\[\spec{A_4D} =\left\{ k^2\pi^2/ [\Re(b_+^{-1/2})]^2
\right\}_{k\in \Z} \subset [0,\infty).\]
\item[c)] If $(a,d) \not \in \{ a^2-ad-1=0 \}$, then
$\spec{A_4D}$ is infinite but it only intersects the real line in
a finite number of points.
\end{itemize}
\end{theorem}
\Proof By virtue of v), $A_4$ is diagonalizable. We assume
$a+d\geq 0$, so that
\[
   \{b_\pm\}\subset \{ \Re(z)\geq 0\}.
\]
Let
\[
  y:=\sqrt{4-(a-d)^2}\qquad  \mathrm{and}\qquad \gam_\pm=a-d \pm
  iy.
\]
Then
\[
   A_4 = \begin{matr2} \gam_+ & \gam_- \\ 2 & 2 \end{matr2}
   \begin{matr2} b_+ & 0 \\ 0 & b_- \end{matr2}
   \begin{matr2} \frac{1}{2iy} & -\frac{\gam_-}{4iy} \\
   -\frac{1}{2iy} & \frac{\gam_+}{4iy} \end{matr2}.
\]
Let $\tet:=\arg{\gam_+}$ and $\alp:=1/\sqrt{b_+}$ so that
$\overline{\alp}=1/\sqrt{b_-}$. Then
\begin{equation} \label{e14}
\begin{aligned}
   \frac{EV(x)}{4\gam_-^2} & =  \frac{2\gam_+}{\gam_-}
   [1-\cos(\alp x)\cos(\overline{\alp}x)]-
        \left(\frac{\gam_+^2}{\gam_-^2} \sqrt{\frac{b_+}{b_-}} +
        \sqrt{\frac{b_-}{b_+}}\right) \sin(\alp x)
        \sin(\overline{\alp}x)  \\
        & =  \frac{2\gam_+}{\overline{\gam_+}}
   [1-\cos(\alp x)\cos(\overline{\alp}x)]-
        \left(\frac{\gam_+^2}{\overline{\gam_+}^2}
        \frac{\overline{\alp}}{\alp}  +
        \frac{\alp}{\overline{\alp}}\right) \sin(\alp x)
        \sin(\overline{\alp}x)  \\
   & =  2\e^{i2\tet}[1-\cos(\alp x)\cos(\overline{\alp}x)]-
        \left( \e^{i4\tet}\frac{\overline{\alp}}{\alp}
         +\frac{\alp}{\overline{\alp}}
        \right) \sin(\alp x)
        \sin(\overline{\alp}x)  \\
   & =  2\e^{i2\tet}[1-\cos(\alp x)\cos(\overline{\alp}x) -
        c \sin(\alp x)\sin(\overline{\alp}x)]  \\
   & =  2\e^{i2\tet} F(x),
\end{aligned}
\end{equation}
where $F(x)$ and
\begin{align*}
   c & :=  \frac{(\overline{\alp}\e^{i2\tet}/\alp)+
   (\alp\e^{-i2\tet}/\overline{\alp})}{2} \\
   & =  \frac{\e^{i(2\tet-2\arg{(\alp))}} +
   \e^{-i(2\tet-2\arg{(\alp)})} }{2} \\
   & =  \cos(2\tet-2\arg(\alp))=\cos(2\tet+\arg{b_+})
\end{align*}
are as in lemma \ref{t16}.

Let us show a). The hypothesis $a-d<0$ ensures $-1< c \leq 1$.
Furthermore $c=1$, if and only if
\[
   \frac{\Im(\gam_+^2)}{\Re(\gam_+^2)}=-\frac{\Im(b_+)}{\Re(b_+)}.
\]
The latter occurs, if and only if
\[
   \frac{y(a-d)}{(a-d)^2-2}
   = -\frac{y}{a+d}.
\]
By simplifying this identity, we gather that $c=1$ for
$a^2-ad-1=0$ which is precisely our assumption. Then, lemma
\ref{t16}-b) and (\ref{e14}) complete the proof of a).

For b), notice that since $a-d>0$, the constant $c$ is now such
that $-1\leq c< 1$ and $c=-1$, if and only if
\[
   \frac{\Im(\gam_+^2)}{\Re(\gam_+^2)}=-\frac{\Im(b_+)}{\Re(b_+)}.
\]
Therefore a similar argument as for a) and lemma \ref{t16}-a) show
this case. In order to prove c) use the fact that $-1<c<1$ in
\[
    R_5 \setminus \{a^2-ad-1=0\},
\]
lemma \ref{t16}-c) and (\ref{e14}). \fin

\gap

\begin{theorem} \label{t29}
In the regions $R_3$ and $R_4$, $\spec{A_4D}$ is infinite, and
\begin{equation*}
\begin{gathered}
   \spec{A_4D} \subset \{ (r+iy_0)^2 : r \in \R \} + [0,\infty) \\
   \spec{A_4D} \subset \{ -(r+iy_0)^2 : r \in \R \} + (-\infty,0]
\end{gathered} \qquad \qquad
  \begin{gathered} (a,d) \in R_3 \\
(a,d) \in R_4, \end{gathered}
\end{equation*}
where in both cases the constant $y_0>0$ only depends upon
$(a,d)$.
\end{theorem}
\Proof We show the result only for $R_3$. According to iii), in
this case $0 < b_- < b_+$ and $A_4$ is diagonalizable. By
expressing the trigonometric functions in exponential form,
\begin{align*}
  EV(x) & = k_1-k_1\cos(\alp x) \cos(\bet x) - k_2 \sin(\alp x)
  \sin (\bet x) \\
  & = k_1+ \frac{k_2-k_1}{4} \left[ \e^{i(\alp+\bet)x} +
  \e^{-i(\alp+\bet)x} \right] + \\
  & \hspace{2in}
  - \frac{k_2+k_1}{4} \left[ \e^{i(\alp-\bet)x} +
  \e^{-i(\alp-\bet)x} \right]
\end{align*}
where $k_1,k_2\in \R$ and $0<\bet<\alp$ are constants we do not
need to specify here. A similar argument involving Hadamard's
theorem as in the proof of lemma \ref{t16} shows that
$\spec{A_4D}$ is infinite.

By putting $x=r+iy$ where $r,y\in \R$, $\gam:=\alp+\bet>0$ and
$\del:=\alp-\bet>0$,
\begin{align*}
   EV(r+iy) & = k_1+ \frac{k_2-k_1}{4}\left[ \e^{-\gam y}\e^{i \gam r} +
  \e^{\gam y}\e^{-i \gam r} \right] + \\ & \hspace{1.5in} - \frac{k_2+k_1}{4}
  \left[ \e^{-\del y}\e^{i \del r} +\e^{\del y}\e^{-i \del r}
  \right].
\end{align*}
Since $\gam>\del>0$, if we chose $y>>0$, the term $\e^{\gam y}$
dominates the expression and so $|EV(r+iy)|\geq c>0$ for a
suitable $c$ independent of $r$. If we chose $y<<0$, the term
$\e^{-\gam y}$ is the one that dominates and again $|EV(r+iy)|$ is
large. This shows that all the zeros of $EV(x)$ must be contained
in a band $\{-y_0 \leq \Im(x) \leq y_0\}$. \fin

The above theorem does not rule out the possibility of negative
eigenvalues when $ad<0$. We will see in the numerical examples,
evidence of points in this region such that $A_4D$ has indeed
negative spectrum.

With regard to finding the minimal $y_0$. We will see in section
\ref{s7} an argument involving Chebyshev polynomial that allows us
to compute in closed form $\spec{A_4D}$ for a certain dense subset
of $R_3$. We will also illustrate this technique in various
examples where the parabolic region is found explicitly.

\gap

If $(a,d)\in R_1$, the matrix $A_4$ is not diagonalizable and so
$EV(x)$ is given by theorem \ref{t26} instead of theorem \ref{t3}.
Nevertheless, similar techniques to the ones we have seen so far
apply to this case.

\begin{lemma} \label{t33} Let $0\not=c \in \R$ and let
\[
   F(x)=x^2 +c [\sin(x)]^2 \qquad \qquad x\in \C.
\]
Then $F(x)$ has an infinite number of zeros in the complex plane
but only a finite number of them lie on $\R$ and on $i\R$.
\end{lemma}
\Proof See the proofs of lemma \ref{t16} and theorem
\ref{t29}.\fin

\begin{theorem} \label{t32} Let $(a,d)\in R_1$. If
$(a,d)=(\pm 1/2,\mp 3/2)$, then \linebreak $\spec{A_4D}=\{0\}$.
Otherwise $\spec{A_4D}$ is infinite but it only intersects the
real line in a finite number of points.
\end{theorem}
\Proof If $a-d=2$,
\[
    A_4=\begin{matr2} b_+ +1 & -1 \\ 1  & b_+-1 \end{matr2} =
    \begin{matr2} 1 & 1 \\ 0 & 1 \end{matr2}
    \begin{matr2} b_+ & 0 \\ 1 & b_+ \end{matr2}
    \begin{matr2} 1&-1 \\ 0 & 1 \end{matr2}
\]
and if $a-d=-2$,
\[
    A_4=\begin{matr2} b_+ -1 & -1 \\ 1  & b_++1 \end{matr2} =
    \begin{matr2} 1 & -1 \\ 0 & 1 \end{matr2}
    \begin{matr2} b_+ & 0 \\ 1 & b_+ \end{matr2}
    \begin{matr2} 1& 1 \\ 0 & 1 \end{matr2}.
\]
Then
\[
   EV(x)=\frac{x^2}{4b_+^3}-\left(1\pm\frac{1}{2b_+}\right)^2\left[\sin\left(
   \frac{x}{\sqrt{b_+}}\right)\right]^2, \qquad \qquad a-d=\pm 2.
\]
The first statement follows from the fact that if $(a,d)=(\pm 1/2,
\mp 3/2)$, then $b_+=\mp 1/2$ and so the trigonometric term
disappear. The second follows from lemma \ref{t33}. \fin

Notice that the curve $a^2-ad-1=0$ meets the region $R_1$ at $(\pm
1/2, \mp 3/2)$. These are the only points where $\spec{A_4D}$ is
finite. Since all self-adjoint operators with compact resolvent
must have an infinite number of eigenvalues, $A_4D$ is not similar
to self-adjoint. All this suggests that for $(a,d)$ in a small
neighbourhood of these points, $\spec{A_4D}$ must be highly
unstable. In the next section we explore closely this idea.


\section{Some numerical results}
\label{s7}

As mentioned previously, this section is devoted to investigating
some aspects of the global spectral evolution of $AD$ when we move
the entries of the matrix $A$. To be more precise, we consider
$A=A_4$ (see section \ref{s6}) and compute $\spec{A_4D}$ as
$(a,d)$ moves along various lines inside $R_1\cup R_3\cup R_5
\subset \R^2$. We also introduce a technique that allows us to
find explicitly $\spec{A_4D}$ when $(a,d)$ are in a certain dense
subset of $R_3$ by computing the roots of certain polynomial
$G(w)$.

Our first task is to decompose $R_3$ into a disjoint union of
curves in order to find the dense subset. For $\alp>1$, let
\[
   d_\pm(a) := \frac{a(\alp^4+1)\pm \sqrt{(\alp^4-1)^2a^2+4\alp^2(\alp^2+1)^2}}
   {2\alp^2}
\]
and let
\[
  \Lam_\pm(\alp) :=\{(a,d_\pm(a)) : a>\mp 1 \}.
\]
Then
\[
    R_3\cap \{a-d<0 \}= \bigcup _{\alp>1} \Lam_+(\alp) \qquad \mathrm{and}
    \qquad R_3 \cap \{a-d>0 \} = \bigcup _{\alp>1} \Lam_- (\alp).
\]
The motivation for this decomposition is found by observing that
for
\[
   A_4 = \begin{matr2} a & -1 \\ 1 & d_\pm(a) \end{matr2},
\]
$\sqrt{b_+/b_-}=\alp$, where $0<b_-<b_+$ are the eigenvalues of
$A_4$. That is, $\Lam_\pm$ are level curves of $\sqrt{b_+/b_-}$ in
the $(a,d)$-plane. Notice that
\[
   R_3= \overline{\bigcup_{1<\alp\in\Q} \Lam_+(\alp) \cup
   \Lam_-(\alp)}.
\]

The key idea behind finding $G(w)$ is that for $(a,d)\in
\Lam_\pm(\alp)$ where $1<\alp\in\Q$, the zeros of the
transcendental function are periodic in the horizontal direction.
We show how to construct this polynomial. The transcendental
function for $A_4D$ is
\begin{align*}
   EV(x)&=k_1\left[1-\cos \left( \frac{x}{\sqrt{b_+}}\right)
   \cos \left(\frac{x}{\sqrt{b_-}}\right) \right] + \\ & \hspace{2in} - k_2
   \sin \left(\frac{x}{\sqrt{b_+}}\right)
   \sin \left(\frac{x}{\sqrt{b_-}}\right)
\end{align*}
where $k_1$ and $k_2$ are two real constants depending upon $a$
and $d$ which we do not need to specify here. Since
\[\sqrt{b_+/b_-} = \alp = p/q, \qquad \qquad p,q\in \Z^+, \]
$\sqrt{b_\pm}$ are rationally related and so the zeros of $EV(x)$
appear periodically in lines parallel to the real axis. By putting
$z=x/(q\sqrt{b_+})$,
\begin{align*}
   EV(x) & = k_1[1-\cos  (pz)\cos (qz) ]  - k_2\sin (pz) \sin (qz)
   \\ &= k_1 + \frac{k_2-k_1}{2} \cos[(p+q)z] - \frac{k_2+k_1}{2}
   \cos[(p-q)z],
\end{align*}
where $p-q<p+q\in \Z^+$. Standard computations show that,
\[
   \cos(m z)=T_m(\cos(z)) \qquad \qquad \qquad m=1,2,\ldots
\]
where $T_m$ a polynomial of degree $m$ (the $m^{th}$ Chebyshev
polynomial of first order). Then by letting
\[
   G(w):= k_1 + \frac{k_2-k_1}{2} T_{(p+q)}(w) - \frac{k_2+k_1}{2}
   T_{(p-q)}(w),
\]
$EV(x)=0$, if and only if $G(\cos(z))=0$. Hence all the zeros of
$EV(x)$ are of the form
\[
   (\pm \arccos(w_0) + 2n\pi)q\sqrt{b_+}\in\C, \qquad \qquad
   n\in \Z
\]
where $w_0$ is a root of $G(w)$. In this manner, $\spec{A_4D}$ is
generated by translations of the roots of $G(w)$.

Although the above method computes  $\spec{A_4D}$ explicitly for
\linebreak $(a,d)\in \Lam_\pm(\alp),1<\alp\in \Q$, its numerical
implementation for large $p+q$ ($>20$ in a PC) is highly unstable
due to the well known instability of the roots of polynomials of
high degree. Nevertheless, no other procedure tried so far, has
proven to be more efficient for estimating large eigenvalues in
$R_3$. Figures \ref{f4}, \ref{f7} and \ref{f8} below were produced
via this approach.

\subsection{Spectral behaviour of $A_4$ for $(a,d)$ close to
$(-1/2,3/2)$} \ \label{ss7.1} \newline By virtue of theorem
\ref{t32}, $\spec{A_4D}=\{0\}$ for $(a,d)=(-1/2,3/2)$. In any
small neighbourhood of this point, the spectrum of $A_4D$ is
infinite so high instability is to be expected. Since $A_4D$ is
holomorphic in $a$ and $d$, every non-zero eigenvalue of $A_4D$
either concentrates at zero or diverges to $\infty$ for $(a,d)\to
(-1/2,3/2)$. We explore this phenomenon in some detail.

According to theorem \ref{t31}-a), if $(a,d) \in R_5$ satisfy
$a^2-ad-1=0$ and $-2< a-d<0$,
\[
  \spec{A_4D}=\left\{ -k^2\pi^2/ [\Im(b_+^{-1/2})]^2\right\}_{k\in \Z },
\]
where $b_+$ as in section \ref{ss6.3}. By taking $a\to -1/2$ and
$d\to 3/2$,
\[
   b_+=\frac{a+d + \sqrt{(a-d)^2-4}}{2}\to 1/2 \in \R,
\]
so that $\Im(b_+^{-1/2})\to 0$. Hence, all non-zero eigenvalues of
$A_4D$ remain negative and escape to $-\infty$ as $(a,d)\in R_5$
approach the critical point on the curve $a^2-ad-1=0$.

\begin{figure}[t]
\begin{picture}(300,300)(50,135) \includegraphics{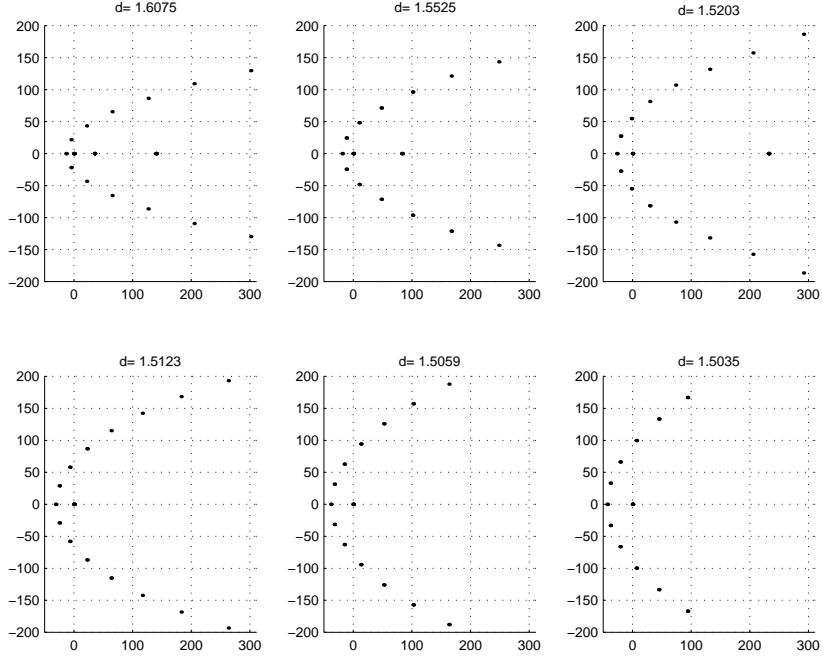}
\end{picture}
\caption{Evolution of the first $16$ eigenvalues of $A_4D$ for
$a=-1/2$ fixed and different values of $d$ close to $3/2$.}
\label{f4}
\end{figure}

In general, not every eigenvalue of $A_4D$ need to be in the left
hand plane when $(a,d)$ is close to $(-1/2,3/2)$. In figure
\ref{f4} we consider the evolution of the first $16$ eigenvalues
of $A_4D$ for $a=-1/2$ fixed and $6$ different values of $d$ from
$d=1.6075$ to $1.5035$. The awkward choice of $d$ correspond to
the sensible values of $\alp\in \Q$; each pair $(-1/2,d)\in
\Lam_+(\alp)$ for $\alp=2,8/5,4/3,5/4,7/6,9/8$. Notice that for
large $p,q$ the polynomial $G(w)$ has $p+q$ roots and nonetheless
all these roots but $0$ lie on the same curve. This curve moves
away from the origin and there is always a negative eigenvalue.
The positive eigenvalues also escape rapidly to $+\infty$ and
there are infinitely many of them.

\begin{figure}[t]
\begin{picture}(300,300)(50,135) \includegraphics{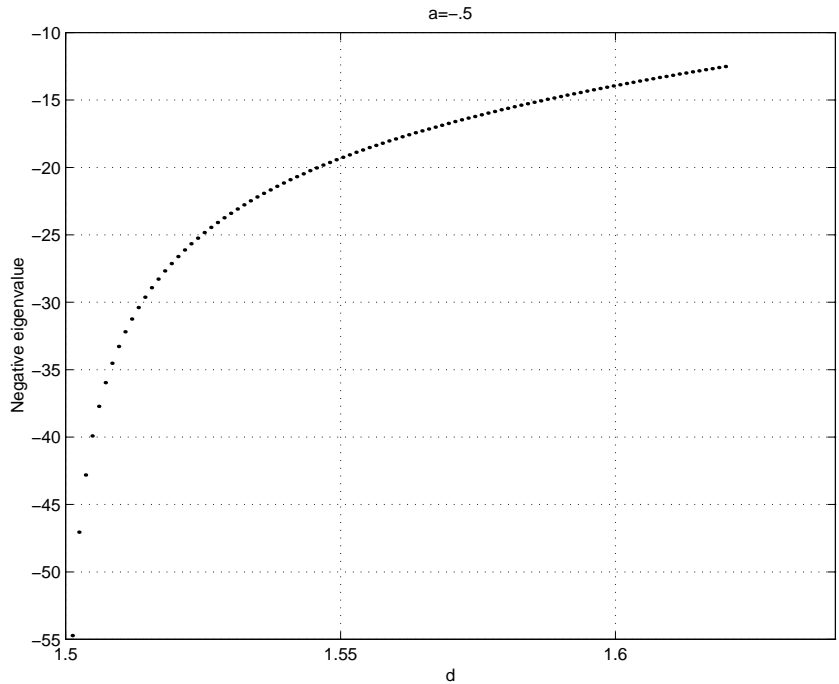}
\end{picture}
\caption{Evolution of the negative eigenvalue for $a=-1/2$ and
$100$ different values of $d$ linearly distributed on the segment
$[1.5012,1.6200]$.} \label{f5}
\end{figure}

In figure \ref{f5} we isolate the negative eigenvalue for $a=-1/2$
against 100 different values of $d$ close to $d=3/2$. This
provides indication of how rapidly it escapes to $-\infty$. In
order to produce this picture, we made use of the algorithm that
Matlab provides to find the zero of $EV(x)$ for $x$ on the
imaginary axis. Comparing with the comment we made earlier in
section \ref{ss6.3}, this provides points in $R_3$ such that
$A_4D$ has a negative eigenvalue of arbitrarily large modulus.

\begin{figure}[t]
\begin{picture}(250,250)(70,130) \includegraphics{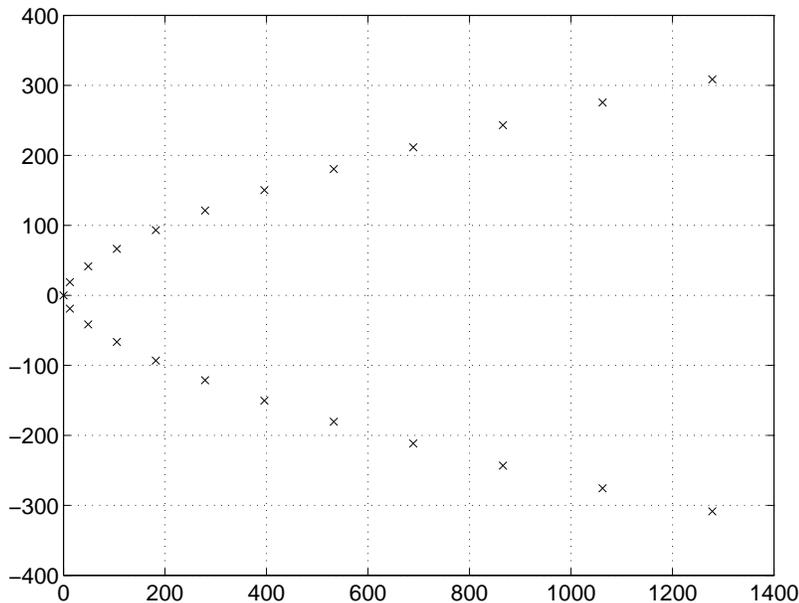}
\end{picture}
\caption{First $23$ eigenvalues of $A_4D$ for $a=0$ and $d=2$.}
\label{f6}
\end{figure}
\subsection{Non-real eigenvalues in $R_1$} \ \label{ss7.2}
\newline We now explore the transition from real to non-real spectrum
by considering the spectral evolution of $A_4$ on the line
\[
   \{(0,d)\in R_3\, :\, d>2\}
\]
close to $(0,2)\in R_1$. In figure \ref{f6} we show the first $23$
eigenvalues of $A_4D$ for $a=0$ and $d=2$. We produced this
graphic by reducing the equation $EV(x)=0$ to a single real
variable and then making use of the algorithm that Maple provides
to find zeros of real functions. According to theorem \ref{t32},
we know that $\spec{A_4D}$ is infinite but there is only finite
intersection with the real line. As the picture suggests, in this
case the origin seems to be the only real eigenvalue.

\begin{figure}[p]
\begin{picture}(300,400)(140,200) \includegraphics{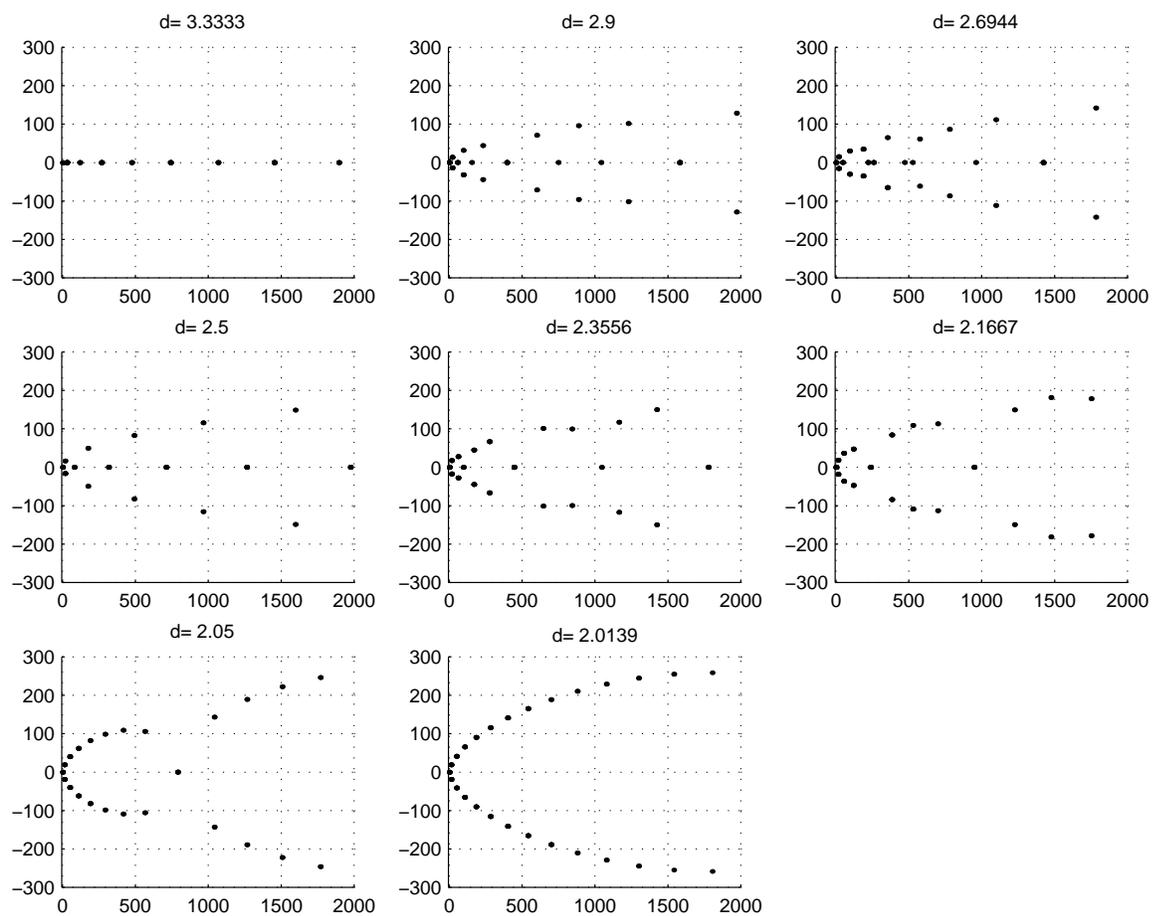}
\end{picture}
\caption{Evolution of the first $17$ eigenvalues (counting
multiplicity) of $A_4D$ for $a=0$ and $d>2$ close to $d=2$.}
\label{f7}
\end{figure}

Figure \ref{f7} shows the evolution of the first $17$ eigenvalues
(counting multiplicity) of $A_4$ when $a=0$ for various different
values of $d$ from $d=3.3333$ to $2.0139$. Each pair $(0,d)\in
\Lam_+(\alp)$ respectively for $\alp=3$, $5/2$, $9/4$, $2$, $9/5$,
$3/2$, $5/4$, $9/8$. The numerical evidence suggests that for
$d=3.3333$ the spectrum is close to the real line and each
eigenvalue is of multiplicity $2$. Each of these operators has
infinitely many real eigenvalues. Unfortunately the method we
employed to find the roots of $G(w)$, is unable to deal with a
finer partition of the $d$-interval. Nonetheless, the global
behaviour of the spectrum can be appreciated, as $d$ approaches to
$2$, each real eigenvalue eventually splits into two conjugate
non-real single eigenvalues stabilizing close to the region in
figure \ref{f6} (see the step $d=2.0139$). Notice that there is no
spectrum in the left hand plane and compare with theorem
\ref{t30}.

\subsection{Spectral evolution close to $R_6$} \ \label{ss7.3}
\newline Another type of peculiar behaviour can be observed as
$(a,d)\in R_3$ approach the region $R_6$, where the matrix $A_4$
is singular and $\spec{A_4D}=\C$. Here we concentrate on the point
$(-1,1) \in R_6$.

\begin{figure}[p]
\begin{picture}(300,400)(140,200) \includegraphics{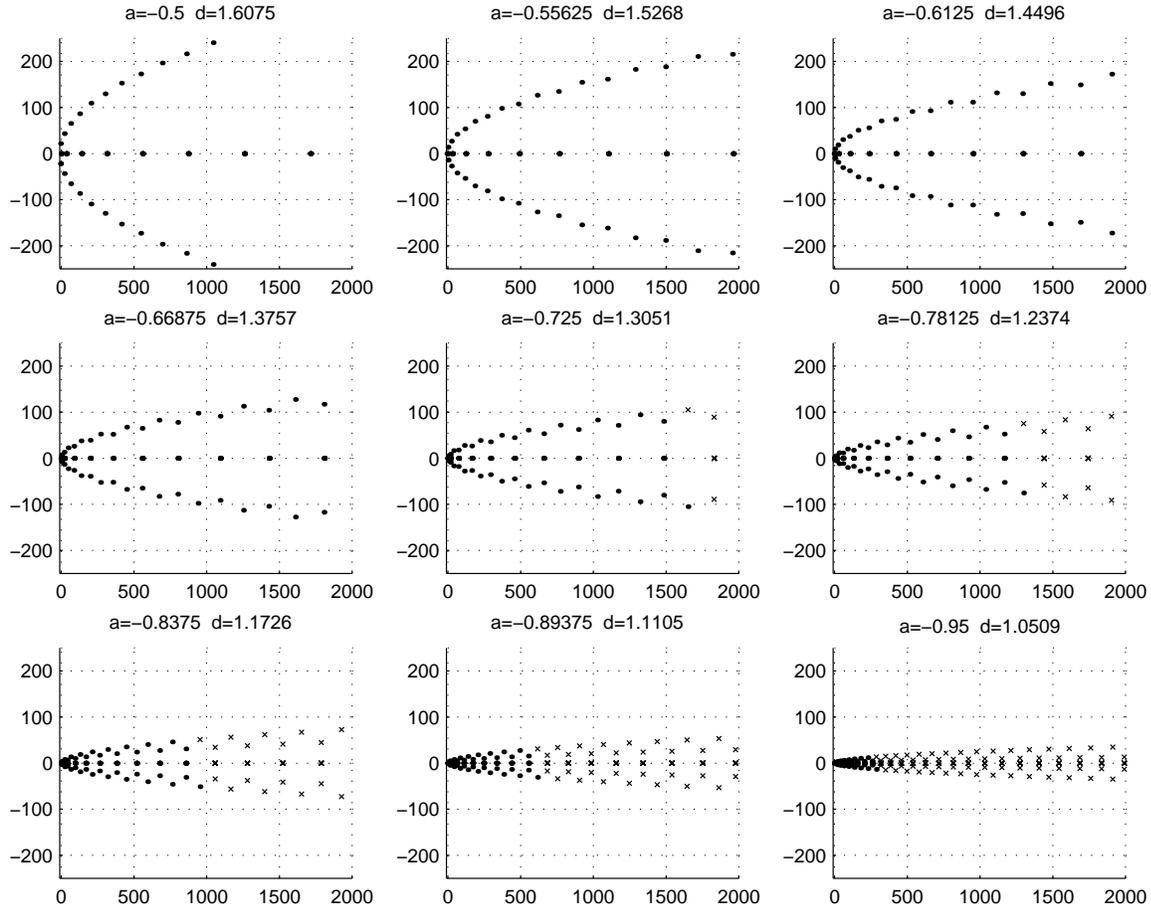}
\end{picture}
\caption{Evolution of the first $100$ eigenvalues of $A_4D$ as
$(a,d)\to (-1,1) \in R_6$ on $\Lam_+(2)$. The dots are the first
$100$ eigenvalues while the crosses the remaining spectrum.}
\label{f8}
\end{figure}

Figure \ref{f8} shows the evolution of the first $100$ eigenvalues
of $A_4D$ (represented by dots) as $(a,d)\in \Lam_+(2)$ approaches
to $(-1,1)\in R_6$. Alongside we also picture the remaining
eigenvalues (represented by crosses) that lie on the box
$[0,2000]\times[-300,300]$ . A very similar behaviour occurs for
$(a,d)\in \Lam_\pm(\alp)$ as $(a,d)\to (\mp 1,\pm 1)\in R_6$ for
other values of $\alp \in \Q$. It can not be appreciated in the
graph provided but there are two conjugate eigenvalues whose real
part is negative. These eigenvalues approach to the origin as
$(a,d)\to (-1,1)$. All the remaining spectrum concentrates on the
real line suggesting that $\spec{A_4D} \to [0,\infty)$ as
$(a,d)\to (-1,1)$ this is in contrast with the fact that
$\spec{AD}=\C$ at $(-1,1)$.

Here we have chosen $p=2$ and $q=1$. This means that $G(w)$ is
only of order $3$ and so the spectrum is always generated by $3$
points. It is not difficult to show analytically that all three
roots converge to $0$ and then rigorously prove that $\spec{A_4D}
\to [0,\infty)$.

\bigskip

{\samepage {\scshape Acknowledgments.} The author wish to thank
Prof. E.~B.~Davies and Prof.~R.~F.~Streater for valuable
discussions of different aspects of this paper. He also would like
to thank Dr.~A.~Aslanyan for many helpful comments and
suggestions.}

\vspace{.3in}

\end{document}